\documentclass[a4paper,11pt,twoside]{article}
\usepackage{amssymb,amsmath,latexsym,geometry,fancyhdr,lineno,hyperref,titletoc}
\contentsmargin{0pt}

\dottedcontents{section}[2em]{\vspace{-1mm}\small}{2em}{0pt}
\dottedcontents{subsection}[5em]{\vspace{-1mm}\small}{3em}{5pt}

\hypersetup{colorlinks,linkcolor= blue,citecolor=blue}
\newtheorem{theorem}{Theorem}[section]

\newtheorem{lemma}[theorem]{Lemma}
\newtheorem{proposition}[theorem]{Proposition}
\newtheorem{remark}{Remark}[section]

\def\proof{\mbox {\textbf{Proof.}~~}}
\numberwithin{equation}{section}
\allowdisplaybreaks[0]
\begin{document}
\title{{\bf\Large Nonradial normalized solutions for nonlinear scalar field equations}}
\author{\\
{ \textbf{\normalsize Louis Jeanjean}}\footnote{For L. Jeanjean this work has been carried out in the framework of the project NONLOCAL (ANR-14-CE25-0013) funded by the French National Research Agency (ANR).}\\
{\it\small Laboratoire de Math\'{e}matiques (CNRS UMR 6623),}\\
{\it\small Universit\'{e} de Bourgogne Franche-Comt\'{e},}\\
{\it\small Besan\c{c}on 25030, France}\\
{\it\small e-mail: louis.jeanjean@univ-fcomte.fr}\\
\\
{ \textbf{\normalsize Sheng-Sen Lu}}\footnote{S.-S. Lu acknowledges the support of the NSF of China (NSFC-11771324), of the China Scholarship Council (CSC-201706250149) and the hospitality of the Laboratoire de Math\'{e}matiques (CNRS UMR 6623), Universit\'{e} de Bourgogne Franche-Comt\'{e}.}\\
{\it\small Center for Applied Mathematics, Tianjin University,}\\
{\it\small Tianjin 300072, PR China}\\
{\it\small e-mail: sslu@tju.edu.cn}}
\date{}
\maketitle
{\bf\normalsize Abstract.} {\small
We study the following nonlinear scalar field equation
  \begin{linenomath*}
    \begin{equation*}
      \left\{
             \begin{aligned}
               -\Delta u&=f(u)-\mu u\quad\text{in}~\mathbb{R}^N,\\
               \|u\|^2_{L^2(\mathbb{R}^N)}&=m,\\
               u&\in H^1(\mathbb{R}^N).
             \end{aligned}
      \right.
    \end{equation*}
  \end{linenomath*}
Here $f\in C(\mathbb{R},\mathbb{R})$, $m>0$ is a given constant and $\mu\in\mathbb{R}$ arises as a Lagrange multiplier. In a mass subcritical case but under general assumptions on the nonlinearity $f$,  we show the existence of one nonradial solution for any $N\geq4$, and obtain multiple (sometimes infinitely many) nonradial solutions when $N=4$ or $N\geq6$. In particular, all these solutions are sign-changing.  }




\pagestyle{fancy}
\fancyhead{} 
\fancyfoot{} 
\renewcommand{\headrulewidth}{0pt}
\renewcommand{\footrulewidth}{0pt}
\fancyhead[CE]{ L. Jeanjean, S.-S. Lu}
\fancyhead[CO]{ Nonradial normalized solutions for nonlinear scalar field equations}
\fancyfoot[C]{\thepage}


\section{Introduction}\label{sect:introduction}
In this paper we consider the nonlinear scalar field equation with an $L^2$ constraint:
  \begin{linenomath*}
    \begin{equation*}\tag{$P_m$}\label{problem}
      \left\{
             \begin{aligned}
               -\Delta u&=f(u)-\mu u\quad\text{in}~\mathbb{R}^N,\\
               \|u\|^2_{L^2(\mathbb{R}^N)}&=m,\\
               u&\in H^1(\mathbb{R}^N).
             \end{aligned}
      \right.
    \end{equation*}
  \end{linenomath*}
Here $N\geq1$, $f\in C(\mathbb{R},\mathbb{R})$, $m>0$ is a given constant and $\mu\in\mathbb{R}$ will arise as a Lagrange multiplier. In particular $\mu\in\mathbb{R}$ does depend on the solution $u \in H^1(\mathbb{R}^N)$ and is not a priori given.

The main feature of \eqref{problem} is that the desired solutions have an a priori prescribed $L^2$-norm. Solutions of this type are often referred to as \emph{normalized solutions}. A strong motivation to study problem \eqref{problem} is that it naturally arises in the search of standing waves of Schr\"{o}dinger type equations of the form
\begin{linenomath*}
    \begin{equation}\label{eq:equation-evolution}
     i \psi_t + \Delta \psi + g(|\psi|^2)\psi =0, \qquad \psi : \mathbb{R}_+ \times \mathbb{R}^N \to \mathbb{C}.
    \end{equation}
  \end{linenomath*}
Here, by standing waves, we mean solutions of \eqref{eq:equation-evolution} of the special form $\psi(t,x) = e^{i\mu t} u(x)$ with $\mu \in \mathbb{R}$ and $u \in H^1(\mathbb{R}^N)$. The study of such type of equations, which had already a strong motivation thirty years ago, see \cite{Be83-1,Lions84-1,Lions84-2}, now lies at the root of several models directly linked with current applications (such as nonlinear optics, the theory of water waves, ...). For these equations, finding solutions with a prescribed $L^2$-norm is particularly relevant since this quantity is preserved along the time evolution. In that direction we refer, in particular, to \cite{C03,CL82,HaSt04,Sh14}. See also the very recent work \cite{St19}.

Under mild assumptions on $f$, it is possible to define the $C^1$ functional
$I: H^1(\mathbb{R}^N)\to\mathbb{R}$ by
  \begin{linenomath*}
    \begin{equation*}\label{eq:functional}
      I(u):=\frac{1}{2}\int_{\mathbb{R}^N}|\nabla u|^2dx-\int_{\mathbb{R}^N}F(u)dx,
    \end{equation*}
  \end{linenomath*}
where $F(t):=\int^t_0f(\tau)d\tau$ for $t\in\mathbb{R}$.  Clearly then, solutions of \eqref{problem} can be characterized as critical points of $I$ submitted to the constraint
  \begin{linenomath*}
    \begin{equation*}\label{eq:constraint}
      S_m:=\left\{u\in H^1(\mathbb{R}^N)~|~\|u\|^2_{L^2(\mathbb{R}^N)}=m\right\}.
    \end{equation*}
  \end{linenomath*}
For future reference, the value $I(u)$ is called the \emph{energy} of $u$.

It is well-known that the study of \eqref{problem} and the type of results one can expect,  do depend on  the behavior of the nonlinearity $f$ at infinity. In particular, this behavior determines whether $I$  is bounded from below on $S_m$. One speaks of a mass subcritical case if $I$ is bounded from below on $S_m$ for any $m>0$, and of a mass supercritical case if $I$ is unbounded from below on $S_m$  for any $m>0$.  One also refers to  a mass critical case when the boundedness from below does depend on the value $m>0$. In this paper we focus on mass subcritical cases and we refer to the papers \cite{BDV13,BS17, Je97} for results in the mass supercritical cases.

The study of the constrained problem \eqref{problem}, or of connected ones, in the mass subcritical case which started with the work of C.A. Stuart in the eighties \cite{St82} had seen a major advance with the work of P.L. Lions \cite{Lions84-1,Lions84-2} on the concentration-compactness principle. Nowadays it is still the object of an intense activity. We refer, in particular, to \cite{CL82, CJS10,HaSt04,JS11,Sh14,St19}.  In these works the authors are mainly interested in the existence of \emph{ground states}, namely of solutions to \eqref{problem} which can be characterized as minimizers of  $I$ among all the solutions. An emphasize is also given to the issue of stability of these solutions, as standing waves of \eqref{eq:equation-evolution}. This is done either, following the strategy laid down in \cite{CL82}, by showing that any minimizing sequences of $I$ on $S_m$ is precompact up to translations \cite{HaSt04,Sh14} or by using more analytic approaches \cite{St19}. Likely, as far as the existence of ground states and their orbital stability is concerned, the most general result is contained in \cite{Sh14}.

Concerning the existence of more than one solution,  the particular case $f(u) = |u|^{\sigma}u$ with $0 < \sigma < 4/N$ and $N \geq 2$ was considered in \cite{Je92} where infinitely many \emph{radial solutions} (with negative energies) were obtained.  For the general result we refer to the recent paper \cite{HT18} by Hirata and Tanaka which still concerns radial solutions. At the end of this paper we shall present the multiplicity result of \cite{HT18} in some details and show that the method we develop in this paper can be used to give an alternative shorter proof of it under a slightly more general setting.

Our aim in the present work is to make further progress in the understanding of the set of solutions to \eqref{problem}.  Roughly speaking, when $N\geq4$, we derive existence and multiplicity results for \emph{nonradial solutions} to \eqref{problem}.
We assume that the nonlinearity $f$ satisfies
\begin{itemize}
    \item[$(f1)$] $f\in C(\mathbb{R},\mathbb{R})$,
    \item[$(f2)$] $\lim_{t\rightarrow0}f(t)/t=0$,
    \item[$(f3)$] $\lim_{t\rightarrow\infty}f(t)/|t|^{q-1}=0$ for some $q < 2^*$ and $\limsup_{t\rightarrow \pm\infty}f(t)t/|t|^{2+4/N}\leq 0$,
    \item[$(f4)$] there exists $\zeta>0$ such that $F(\zeta)>0$,
    \item[$(f5)$] $f(-t)=-f(t)$ for all $t\in\mathbb{R}$.
  \end{itemize}
We shall also make use of the following condition
  \begin{linenomath*}
    \begin{equation}\label{eq:f_key1}
      \lim_{t\to0}\frac{F(t)}{|t|^{2+\frac{4}{N}}}=+\infty,
    \end{equation}
\end{linenomath*}
which is originally introduced in \cite{Lions84-2}; see also \cite{Sh14}. As a simple example of the nonlinearity satisfying $(f1)-(f5)$ (and also \eqref{eq:f_key1}) we have
  \begin{linenomath}
    \begin{equation*}
      f(t)=|t|^{p-2}t-|t|^{q-2}t\qquad\text{with}~ 2< p <2+\frac{4}{N} < q <2^*.
    \end{equation*}
  \end{linenomath}
								
To state our results, we introduce some notations. Assume that $N\geq4$ and $2\leq M\leq N/2$. Let us fix $\tau\in \mathcal{O}(N)$ such that $\tau(x_1,x_2,x_3)=(x_2,x_1,x_3)$ for $x_1,x_2\in\mathbb{R}^M$ and $x_3\in\mathbb{R}^{N-2M}$, where $x=(x_1,x_2,x_3)\in\mathbb{R}^N=\mathbb{R}^M\times\mathbb{R}^M\times\mathbb{R}^{N-2M}$. We define
  \begin{linenomath*}
    \begin{equation*}
      X_\tau:=\left\{u\in H^1(\mathbb{R}^N)~|~u(\tau x)=-u(x)~\text{for all}~x\in\mathbb{R}^N\right\}.
    \end{equation*}
  \end{linenomath*}
It is clear that $X_\tau$ does not contain nontrivial radial functions. Let $H^1_{\mathcal{O}_1}(\mathbb{R}^N)$ denote the subspace of invariant functions with respect to $\mathcal{O}_1$, where $\mathcal{O}_1:=\mathcal{O}(M)\times\mathcal{O}(M)\times \text{id}\subset \mathcal{O}(N)$ acts isometrically on $H^1(\mathbb{R}^N)$. We also consider $\mathcal{O}_2:=\mathcal{O}(M)\times\mathcal{O}(M)\times\mathcal{O}(N-2M)\subset \mathcal{O}(N)$ acting isometrically on $H^1(\mathbb{R}^N)$ with the subspace of invariant functions denoted by $H^1_{\mathcal{O}_2}(\mathbb{R}^N)$. Here we agree that the components corresponding to $N-2M$ do not exist when $N=2M$. Clearly, $H^1_{\mathcal{O}_2}(\mathbb{R}^N)$ is in general a subspace of $H^1_{\mathcal{O}_1}(\mathbb{R}^N)$, but $H^1_{\mathcal{O}_2}(\mathbb{R}^N)= H^1_{\mathcal{O}_1}(\mathbb{R}^N)$ when $N=2M$.
\smallskip

For notational convenience, we set $X_1:= H^1_{\mathcal{O}_1}(\mathbb{R}^N)\cap X_\tau$ and $X_2:= H^1_{\mathcal{O}_2}(\mathbb{R}^N)\cap X_\tau$. Our first main result concerns the existence of one nonradial solution to \eqref{problem}.
  \begin{theorem}\label{theorem:nonradialsolution}
    Assume that $N\geq4$ and $f$ satisfies $(f1)-(f5)$. Define
      \begin{linenomath*}
        \begin{equation*}\label{eq:infimum1}
          E_m:=\inf_{u\in S_m\cap X_1}I(u).
        \end{equation*}
      \end{linenomath*}
    Then $E_m > - \infty$ and the mapping $m \mapsto E_m$ is nonincreasing and continuous. Moreover
		\begin{itemize}
        \item[$(i)$] there exists a uniquely determined number $m^*\in[0,\infty)$ such that
      \begin{linenomath*}
        \begin{equation*}
          E_m=0\quad\text{if}~0<m\leq m^*,\qquad E_m<0\quad\text{when}~m>m^*;
        \end{equation*}
      \end{linenomath*}
        \item[$(ii)$] when $m>m^*$, the infimum $E_m$ is reached and thus \eqref{problem} has one nonradial solution $w\in X_1$ such that $I(w)=E_m$;
				\item[$(iii)$] when $0 < m <m^*$, $E_m$ is not reached;
				\item[$(iv)$] $m^*=0$ if in addition \eqref{eq:f_key1} holds.
      \end{itemize}
  \end{theorem}

Our second main result concerns the multiplicity of nonradial solutions to \eqref{problem}. Let $\Sigma(S_m\cap X_2)$ be the family of closed symmetric subsets of $S_m\cap X_2$,  that is
  \begin{linenomath*}
    \begin{equation*}
      \Sigma(S_m\cap X_2):= \big\{ A \subset S_m\cap X_2 ~|~A~\text{is closed},~-A= A\big\},
    \end{equation*}
  \end{linenomath*}
and denote by $\mathcal{G}(A)$ the genus of $A\in \Sigma(S_m\cap X_2)$. For the definition of the genus and its basic properties, one may refer to Section \ref{sect:minimaxtheorem}.
  \begin{theorem}\label{theorem:nonradialsolutions}
    Assume that $N\geq4$, $N-2M\neq1$ and $f$ satisfies $(f1)-(f5)$. Define the minimax values
      \begin{linenomath*}
        \begin{equation*}
          E_{m,k}:=\inf_{A\in \Gamma_{m,k}}\sup_{u\in A}I(u),
        \end{equation*}
      \end{linenomath*}
   where $\Gamma_{m,k}:=\{A\in \Sigma(S_m\cap X_2)~|~\mathcal{G}(A)\geq k\}$.  Then the following statements hold.
      \begin{itemize}
		\item[$(i)$] $-\infty<E_{m,k}\leq E_{m,k+1}\leq0$ for all $m>0$ and $k\in\mathbb{N}$.
		\item[$(ii)$] For any $k\in\mathbb{N}$, the mapping $m\mapsto E_{m,k}$ is nonincreasing and continuous.
        \item[$(iii)$] For each $k\in\mathbb{N}$, there exists a uniquely determined $m_k\in[0,\infty)$ such that
				\begin{linenomath*}
        \begin{equation*}
          E_{m,k}=0\quad\text{if}~0<m\leq m_k,\qquad E_{m,k}<0\quad\text{when}~m>m_k.
        \end{equation*}
      \end{linenomath*}
			When $m>m_k$,
                      \eqref{problem} has $k$ distinct nonradial solutions belonging to $X_2$ and associated to the levels $E_{m,j}$ ($j=1,2, \cdots, k$).
        \item[$(iv)$] Assume in addition \eqref{eq:f_key1}, then $m_k =0$ for any $k \in \mathbb{N}$ and thus  \eqref{problem} has infinitely many nonradial solutions $\{w_k\}^\infty_{k=1}\subset X_2$ for all $m>0$. In particular, $I(w_k)=E_{m,k}<0$ for each $k\in\mathbb{N}$ and $I(w_k)\to0$ as $k\to\infty$.
      \end{itemize}
  \end{theorem}

The question of the existence of nonradial solutions to the free equation
\begin{linenomath*}
    \begin{equation}\label{eq:equation-free}
      - \Delta u = f(u) - \mu u, \qquad u \in H^1(\mathbb{R}^N)
    \end{equation}
  \end{linenomath*}
was raised in \cite[Section 10.8]{Be83-2} and remained open  for a long time. Partial results, namely for specific nonlinearities $f$, were first obtained by Bartsch and Willem \cite{Ba93}  (we refer to \cite{Lions86} if nonradial complex solutions are of interest to the reader).  The authors worked in dimension $N=4$ and $N \geq 6$ assuming an Ambrosetti-Rabinowitz type condition. Actually the idea of considering subspaces of $H^1(\mathbb{R}^N)$ as $H^1_{\mathcal{O}_2}(\mathbb{R}^N)$ originates from \cite{Ba93}. Note also the work \cite{Lo04} in which the problem is solved when $N=5$ by introducing the $\mathcal{O}_1$ action on $H^1(\mathbb{R}^N)$. Finally it was only very recently that, under general assumptions on $f$, a positive answer to the existence and multiplicity of nonradial solutions was given \cite{Me17}. Note also that in \cite{JL18} the authors gave an alternative proof of the results of \cite{Me17} with more elementary arguments. All these results however consider the equation \eqref{eq:equation-free} without prescribing the $L^2$-norm of the solutions.

The present paper is, up to our knowledge, the first to consider the existence of nonradial normalized solutions.  We also observe that the nonradial solutions given by Theorems \ref{theorem:nonradialsolution} and \ref{theorem:nonradialsolutions} change signs. In sharp constrast to the unconstrained case \eqref{eq:equation-free} where numerous results have been established, see for example \cite{Ba93_1,LW08,MPW12},  the existence of \emph{sign-changing solutions} had not been studied yet  for $L^2$-constrained problems.

Let us now give some ideas of the proofs of Theorems  \ref{theorem:nonradialsolution} and \ref{theorem:nonradialsolutions}.  To prove the multiplicity result stated in Theorem \ref{theorem:nonradialsolutions}, we work in the space $X_2:= H^1_{\mathcal{O}_2}(\mathbb{R}^N)\cap X_\tau$ and make use of classical minimax arguments (see Theorem \ref{theorem:minimax} below). Since $N \geq 4$ and $N-2M \neq 1$, we can benefit from the compact inclusion $X_2\hookrightarrow L^p(\mathbb{R}^N)$  for all $2<p<2N/(N-2)$. This result, which is due to P. L. Lions \cite{Lions82}, allows to show that $I_{|S_m\cap X_2}$ satisfies the Palais-Smale condition at any level $c<0$, see Lemma \ref{lemma:PS}. We observe that the proof of Lemma \ref{lemma:PS}, and likely its conclusion, fail at level $c\geq 0$.  Thus another key point is to verify that the minimax levels $E_{m,k}$ are indeed negative for some or any $m>0$ and $k\in\mathbb{N}$.  Relying on the construction of some special mappings done in \cite{Be83-2,JL18, Me17}, we manage to do this in Lemmas \ref{lemma:geo2} and \ref{lemma:Emk}.  Note also that to derive the existence of $m_k$ in Theorem \ref{theorem:nonradialsolutions} $(iii)$, we need Theorem \ref{theorem:nonradialsolutions} $(ii)$ which is proved in Lemma \ref{lemma:Emk}.

For the proof of Theorem \ref{theorem:nonradialsolution}, we work in the space $X_1:= H^1_{\mathcal{O}_1}(\mathbb{R}^N)\cap X_\tau$. In the case where $N-2M=0$, we have $X_1=X_2$ (with $N-2M\neq1$). Since then  $E_m= E_{m,1}$ and $m^* = m_1$, the existence of a minimizer for $E_m$ when $m>m^*$  follows directly from Theorem \ref{theorem:nonradialsolutions}.  When $N-2M\neq0$, the inclusion $X_1\hookrightarrow L^p(\mathbb{R}^N)$ is not compact for any $2<p<2N/(N-2)$ and the Palais-Smale condition does not hold any more. To derive our existence result in this case, using concentration compactness type arguments, we study carefully the behavior of the minimizing sequences of $E_m$. Here again it is essential to know in advance that the suspected critical level is negative.

\begin{remark}\label{remark:stability}
If the stability of the ground states, as studied for example in \cite{HaSt04, Sh14, St19}, is by now relatively well understood, the issue of the orbital stability, or more likely orbital instability, of the other critical points of $I$ restricted to $S_m$ is still totally open. We believe an interesting but challenging question would be to prove that the solution obtained in Theorem \ref{theorem:nonradialsolution}, which enjoys a well defined variational characterization, is orbitally unstable.
\end{remark}

\begin{remark}\label{remark:extension}
Taking advantage of an idea first introduced in \cite{Be83-1}, it is possible to find solutions to \eqref{problem} under $(f1)-(f5)$ when $(f3)$ is replaced by the more general condition
\begin{itemize}
		\item[$(f3)'$] $\limsup_{t\rightarrow +\infty}f(t)/t^{1+4/N} \leq 0$.
  \end{itemize}
Indeed, assume that $f$ satisfies $(f1)$, $(f2)$, $(f3)'$, $(f4)$ and $(f5)$. If $f(t)\geq0$ for all $t\geq\zeta$, then $f$ satisfies $(f1)-(f5)$. Otherwise, we set
\begin{linenomath*}
\begin{equation*}
\zeta_1:=\inf\{t\geq\zeta~|~f(t)=0\}\qquad\text{and}\qquad
\widetilde{f}(t):=\left\{
\begin{aligned}
&f(t),~&\text{for}~|t|\leq \zeta_1,\\
&0,&\text{for}~|t|>\zeta_1.\\
\end{aligned}
\right.
\end{equation*}
\end{linenomath*}
Clearly $\widetilde{f}$ satisfies $(f1)-(f5)$.  Also, for any couple $(u, \mu)\in S_m \times \mathbb{R}_+$ satisfying
\begin{linenomath*}
  \begin{equation*}
    -\Delta{u}= \widetilde{f}(u)-\mu u\quad\text{in}~\mathbb{R}^N,
  \end{equation*}
\end{linenomath*}
 the strong maximum principle tells us that $|u(x)|\leq \zeta_1$ for all $x\in\mathbb{R}^N$ and so $u \in S_m$ actually satisfies $- \Delta u = f(u) - \mu u$ in $\mathbb{R}^N$. Applying Theorems  \ref{theorem:nonradialsolution} and \ref{theorem:nonradialsolutions} with $\widetilde{f}$ and noting that the Lagrange multipliers associated to the solutions obtained by these theorems belong to $\mathbb{R}_+$ (see the proof of Lemma \ref{lemma:PS}), we thus obtain existence and multiplicity results to \eqref{problem}.  Note however that under $(f3)'$ the functional $I$ is in general not more defined and so there is no direct connection between our solutions and the evolution equation \eqref{eq:equation-evolution}.
 \end{remark}

\begin{remark}\label{remark:partially_radial}
  Since we work in the spaces $X_{1}$ and $X_{2}$, the nonradial solutions we obtain are still partially radial, that is, radial with respect to certain groups of directions. Actually, in order to keep a minimal compactness, we find necessary to impose that the subspaces in which radial symmetry is preserved are at least two dimensional. This implies the condition $M \geq 2$ and in turn that $N \geq 4$. In addition, in the definition of $X_1$ and $X_2$, we introduce the parity property $X_{\tau}$ with respect to a certain ``diagonal". It is this odd property which ensures that the solutions are not globally radial. As an open problem, it would be interesting to inquire if there exist nonradial solutions of \eqref{problem} having less symmetry or directly living in $\mathbb{R}^2$ or $\mathbb{R}^3$.
\end{remark}

The paper is organized as follows. In Section \ref{sect:minimaxtheorem} we present the version of the minimax theorem that will be used in the proof of Theorem \ref{theorem:nonradialsolutions}.  Section \ref{sect:preliminaries} establishes  some key technical  points to be used in the proofs of the main results. In Section \ref{sect:proofs} we prove  Theorems \ref{theorem:nonradialsolution} and \ref{theorem:nonradialsolutions}. Finally, with the approach used to prove Theorem \ref{theorem:nonradialsolutions}, we prove in Section \ref{sect:theoremB} a slight extension of the multiplicity result due to Hirata and Tanaka \cite[Theorem 0.2]{HT18}.

\section{A minimax theorem}\label{sect:minimaxtheorem}
In this section, we present a minimax theorem for a class of constrained even functionals. Let us point out that closely related results do exist in the literature, see in particular e.g., \cite[Section 8]{Be83-2}, \cite{Ra86}, \cite{Sz88} and \cite[Chapter 5]{Wi96}. The present version is well suited to deal with the nonlinear scalar field equations considered in this paper.

To formulate the minimax theorem, we need some notations. Let $\mathcal{E}$ be a real Banach space with norm $\|\cdot\|_\mathcal{E}$ and $\mathcal{H}$ be a real Hilbert space with inner product $(\cdot,\cdot)_\mathcal{H}$. We identify  $\mathcal{H}$ with its dual space and assume that $\mathcal{E}$ is embedded continuously in $\mathcal{H}$. For any $m>0$, define the manifold
  \begin{linenomath*}
    \begin{equation*}
      \mathcal{M}:=\{u\in \mathcal{E}~|~(u,u)_\mathcal{H}=m\},
    \end{equation*}
  \end{linenomath*}
which is endowed with the topology inherited from $\mathcal{E}$. Clearly, the tangent space of $\mathcal{M}$ at a point $u\in\mathcal{M}$ is defined by
  \begin{linenomath*}
    \begin{equation*}
      T_u\mathcal{M}:=\{v\in \mathcal{E}~|~(u,v)_\mathcal{H}=0\}.
    \end{equation*}
  \end{linenomath*}
Let $I\in C^1(\mathcal{E},\mathbb{R})$, then $I_{|\mathcal{M}}$ is a functional of class $C^1$ on $\mathcal{M}$. The norm of the derivative of $I_{|\mathcal{M}}$ at any point $u\in\mathcal{M}$ is defined by
  \begin{linenomath*}
    \begin{equation*}
      \|I_{|\mathcal{M}}'(u)\|:=\sup_{\|v\|_\mathcal{E}\leq 1,~v\in T_u\mathcal{M}}|\langle I'(u),v\rangle|.
    \end{equation*}
  \end{linenomath*}
A point $u\in\mathcal{M}$ is said to be a critical point of $I_{|\mathcal{M}}$ if $I'_{|\mathcal{M}}(u)=0$ (or, equivalently, $\|I_{|\mathcal{M}}'(u)\|=0$). A number $c\in\mathbb{R}$ is called a critical value of $I_{|\mathcal{M}}$ if $I_{|\mathcal{M}}$ has a critical point $u\in\mathcal{M}$ such that $c=I(u)$. We say that $I_{|\mathcal{M}}$ satisfies the Palais-Smale condition at a level $c\in\mathbb{R}$, $(PS)_c$ for short, if any sequence $\{u_n\}\subset \mathcal{M}$ with $I(u_n)\to c$ and $\|I'_{|\mathcal{M}}(u_n)\|\to 0$ contains a convergent subsequence.

Noting that $\mathcal{M}$ is symmetric with respect to $0\in\mathcal{E}$ and $0\notin\mathcal{M}$, we introduce the notation of the genus. Let $\Sigma(\mathcal{M})$ be the family of closed symmetric subsets of $\mathcal{M}$. For any nonempty set $A\in \Sigma(\mathcal{M})$, the genus $\mathcal{G}(A)$ of $A$ is defined as the least integer $k\geq1$ for which there exists an odd continuous mapping $\varphi:A\to\mathbb{R}^k\setminus\{0\}$. We set $\mathcal{G}(A)=\infty$ if such an integer does not exist, and set $\mathcal{G}(A)=0$ if $A=\emptyset$.  For each $k\in\mathbb{N}$, let $\Gamma_k:=\{A\in \Sigma(\mathcal{M})~|~\mathcal{G}(A)\geq k\}$.

We now state the minimax theorem and then give a detailed proof for completeness.
  \begin{theorem}[Minimax theorem]\label{theorem:minimax}
    Let $I:\mathcal{E}\to\mathbb{R}$ be an even functional of class $C^1$. Assume that $I_{|\mathcal{M}}$ is bounded from below and satisfies the $(PS)_c$ condition for all $c<0$,  and that $\Gamma_k \neq \emptyset $ for each $k \in \mathbb{N}$.   Then a sequence of minimax values $-\infty<c_1\leq c_2\leq\cdots\leq c_k\leq \cdots$ can be defined as follows:
      \begin{linenomath*}
        \begin{equation*}
          c_k:=\inf_{A\in\Gamma_k}\sup_{u\in A}I(u),\qquad k\geq1,
        \end{equation*}
      \end{linenomath*}
    and the following statements hold.
      \begin{itemize}
        \item[$(i)$] $c_k$ is a critical value of $I_{|\mathcal{M}}$ provided $c_k<0$.
        \item[$(ii)$] Denote by $K^c$ the set of critical points of $I_{|\mathcal{M}}$ at a level $c\in\mathbb{R}$. If
              \begin{linenomath*}
                \begin{equation*}
                  c_k=c_{k+1}=\cdots=c_{k+l-1}=:c<0\qquad\text{for some}~k,l\geq1,
                \end{equation*}
              \end{linenomath*}
            then $\mathcal{G}(K^c)\geq l$. In particular, $I_{|\mathcal{M}}$ has infinitely many critical points at the level $c$ if $l\geq2$.
        \item[$(iii)$] If $c_k<0$ for all $k\geq1$, then $c_k\to0^-$ as $k\to\infty$.
      \end{itemize}
  \end{theorem}

To prove Theorem \ref{theorem:minimax}, we shall need some basic properties of the genus. For $A\subset\mathcal{M}$ and $\delta>0$, denote by $A_\delta$ the uniform $\delta$-neighborhood of $A$ in $\mathcal{M}$, that is,
  \begin{linenomath*}
    \begin{equation*}
      A_\delta:=\{u\in\mathcal{M}~|~\inf_{v\in A}\|u-v\|_\mathcal{E}\leq\delta\}.
    \end{equation*}
  \end{linenomath*}
Since $\mathcal{M}$ is a closed symmetric subset of $\mathcal{E}$, repeating the arguments in \cite[Section 7]{Ra86}, one can get Proposition \ref{proposition:genus} below which is sufficient for our use.
  \begin{proposition}\label{proposition:genus}
    Let $A,B\in\Sigma(\mathcal{M})$. Then the following statements hold.
      \begin{itemize}
        \item[$(i)$] If $\mathcal{G}(A)\geq2$, then $A$ contains infinitely many distinct points.
        \item[$(ii)$] $\mathcal{G}(\overline{A\setminus B})\geq \mathcal{G}(A)-\mathcal{G}(B)$ if $\mathcal{G}(B)<\infty$.
        \item[$(iii)$] If there exists an odd continuous mapping $\psi:\mathbb{S}^{k-1}\to A$, then $\mathcal{G}(A)\geq k$.
        \item[$(iv)$] If $A$ is compact, then $\mathcal{G}(A)<\infty$ and there exists $\delta>0$ such that $A_\delta\in\Sigma(\mathcal{M})$ and $\mathcal{G}(A_\delta)=\mathcal{G}(A)$.
      \end{itemize}
  \end{proposition}

We shall also need the following quantitative deformation lemma whose proof is similar to that of \cite[Lemma 2.3]{Wi96}. For $c<d$, set $I^c_{|\mathcal{M}}:=\{u\in\mathcal{M}~|~I(u)\leq c\}$ and $I^{-1}_{|\mathcal{M}}([c,d]):=\{u\in\mathcal{M}~|~c\leq I(u)\leq d\}$.
  \begin{lemma}\label{lemma:QDL}
    Assume $I_{|\mathcal{M}}\in C^1(\mathcal{M},\mathbb{R})$. Let $S\subset\mathcal{M}$, $c\in\mathbb{R}$, $\varepsilon>0$ and $\delta>0$ such that
      \begin{linenomath*}
        \begin{equation}\label{eq:QDL}
          \|I'_{|\mathcal{M}}(u)\|\geq \frac{8\varepsilon}{\delta}\qquad\text{for all}~u\in I^{-1}_{|\mathcal{M}}([c-2\varepsilon,c+2\varepsilon])\cap S_{2\delta}.
        \end{equation}
      \end{linenomath*}
    Then there exists a mapping $\eta\in C([0,1]\times \mathcal{M},\mathcal{M})$ such that
      \begin{itemize}
        \item[$(i)$] $\eta(t,u)=u$ if $t=0$ or if $u\not\in I^{-1}_{|\mathcal{M}}([c-2\varepsilon,c+2\varepsilon])\cap S_{2\delta}$,
        \item[$(ii)$] $\eta(1,I^{c+\varepsilon}_{|\mathcal{M}}\cap S)\subset I^{c-\varepsilon}_{|\mathcal{M}}$,
        \item[$(iii)$] $I(\eta(t,u))$ is nonincreasing in $t\in[0,1]$ for any $u\in\mathcal{M}$,
        \item[$(iv)$] $\eta(t,u)$ is odd in $u\in\mathcal{M}$ for any $t\in[0,1]$ if $I_{|\mathcal{M}}$ is even,
        \item[$(v)$] $\eta(t, \cdot)$ is a homeomorphism of $\mathcal{M}$ for each $t \in [0,1]$.
      \end{itemize}
  \end{lemma}

With Proposition \ref{proposition:genus} and Lemma \ref{lemma:QDL} in hand, we can now prove Theorem \ref{theorem:minimax}.

\medskip
\noindent
\textbf{Proof of Theorem \ref{theorem:minimax}.} Item $(i)$ is a special case of Item $(ii)$ when $l=1$, so we go straight to the proof of Item $(ii)$. Obviously, $K^c\in \Sigma(\mathcal{M})$ and $K^c$ is compact by the $(PS)_c$ condition. If $\mathcal{G}(K^c)\leq l-1$, by Proposition \ref{proposition:genus} $(iv)$, there exists $\delta>0$ such that
  \begin{linenomath*}
    \begin{equation*}
      \mathcal{G}(K^c_{3\delta})=\mathcal{G}(K^c)\leq l-1.
    \end{equation*}
  \end{linenomath*}
We remark here that $K^c_{3\delta}=\emptyset$ if $K^c=\emptyset$. Let $S:=\overline{\mathcal{M}\setminus K^c_{3\delta}}\subset \mathcal{M}$. Clearly, there exists $\varepsilon>0$ small enough such that \eqref{eq:QDL} holds (if not, one will get a Palais-Smale sequence $\{u_n\}\subset S_{2\delta}$ of $I_{|\mathcal{M}}$ at the level $c<0$ and thus an element $v\in S_{2\delta}\cap K^c$ by the $(PS)_c$ condition, which leads a contradiction since $S_{2\delta}\cap K^c=\emptyset$). Therefore, Lemma \ref{lemma:QDL} yields a mapping $\eta\in C([0,1]\times\mathcal{M},\mathcal{M})$ such that
  \begin{linenomath*}
    \begin{equation*}
      \eta(1,I^{c+\varepsilon}_{|\mathcal{M}}\cap S)\subset I^{c-\varepsilon}_{|\mathcal{M}}\qquad\text{and}\qquad \eta(t,\cdot)~\text{is odd for all}~t\in[0,1].
    \end{equation*}
  \end{linenomath*}
Choose $A\in\Gamma_{k+l-1}$ such that $\sup_{u\in A}I(u)\leq c+\varepsilon$. It is clear that $\overline{A\setminus K^c_{3\delta}}\subset I^{c+\varepsilon}_{|\mathcal{M}}\cap S$ and thus
  \begin{linenomath*}
    \begin{equation}\label{eq:minimax}
      \eta(1,\overline{A\setminus K^c_{3\delta}})\subset\eta(1,I^{c+\varepsilon}_{|\mathcal{M}}\cap S)\subset I^{c-\varepsilon}_{|\mathcal{M}}.
    \end{equation}
  \end{linenomath*}
On the other hand, since $\mathcal{G}(\overline{A\setminus K^c_{3\delta}})\geq \mathcal{G}(A)-\mathcal{G}(K^c_{3\delta})\geq k$ by Proposition \ref{proposition:genus} $(ii)$, we have $\overline{A\setminus K^c_{3\delta}}\in\Gamma_k$ and then $\eta(1,\overline{A\setminus K^c_{3\delta}})\in\Gamma_k$. Now, by the definition of $c_k$ and \eqref{eq:minimax}, we get a contradiction:
  \begin{linenomath*}
    \begin{equation*}
      c=c_k\leq \sup_{u\in\eta(1,\overline{A\setminus K^c_{3\delta}})}I(u)\leq c-\varepsilon.
    \end{equation*}
  \end{linenomath*}
Thus $\mathcal{G}(K^c)\geq l$. In view of Proposition \ref{proposition:genus} $(i)$, we complete the proof of Item $(ii)$.

To prove Item $(iii)$, we assume by contradiction that there exists $c<0$ such that $c_k\leq c$ for all $k\geq 1$ and $c_k\to c$ as $k\to\infty$. By the $(PS)_c$ condition,  $K^c$ is a (symmetric) compact set.  Thus, by Proposition \ref{proposition:genus} $(iv)$, there exists $\delta>0$ such that
  \begin{linenomath*}
    \begin{equation*}
      \mathcal{G}(K^c_{3\delta})=\mathcal{G}(K^c)=:q<\infty.
    \end{equation*}
  \end{linenomath*}
Let $S:=\overline{\mathcal{M}\setminus K^c_{3\delta}}\subset\mathcal{M}$. Since \eqref{eq:QDL} holds for small enough $\varepsilon>0$, we know from Lemma \ref{lemma:QDL} that a mapping $\eta\in C([0,1]\times\mathcal{M},\mathcal{M})$ exists such that $\eta(1,I^{c+\varepsilon}_{|\mathcal{M}}\cap S)\subset I^{c-\varepsilon}_{|\mathcal{M}}$ and $\eta(t,\cdot)$ is odd for any $t\in[0,1]$. Choose $k\geq1$ large enough such that $c_k>c-\varepsilon$ and take $A\in\Gamma_{k+q}$ such that $\sup_{u\in A}I(u)\leq c_{k+q}+\varepsilon$. Noting that $c_{k+q}\leq c$, we have $\overline{A\setminus K^c_{3\delta}}\subset I^{c+\varepsilon}_{|\mathcal{M}}\cap S$ and thus
  \begin{linenomath*}
    \begin{equation*}
      \eta(1,\overline{A\setminus K^c_{3\delta}})\subset\eta(1,I^{c+\varepsilon}_{|\mathcal{M}}\cap S)\subset I^{c-\varepsilon}_{|\mathcal{M}}.
    \end{equation*}
  \end{linenomath*}
On the other hand, since $\mathcal{G}(\overline{A\setminus K^c_{3\delta}})\geq \mathcal{G}(A)-\mathcal{G}(K^c_{3\delta})\geq k$, we have $\overline{A\setminus K^c_{3\delta}}\in\Gamma_k$ and then $\eta(1,\overline{A\setminus K^c_{3\delta}})\in\Gamma_k$. We now reach a contradiction:
  \begin{linenomath*}
    \begin{equation*}
      c_k\leq \sup_{u\in\eta(1,\overline{A\setminus K^c_{3\delta}})}I(u)\leq c-\varepsilon,
    \end{equation*}
  \end{linenomath*}
for we chosen $k$ large enough such that $c_k>c-\varepsilon$. Thus $c_k\to0^-$ as $k\to\infty$.~~$\square$

To end this section, we recall a characterization result in \cite{Be83-2} which allows to check the $(PS)_c$ condition in a convenient way.
  \begin{lemma}[{\cite[Lemma 3]{Be83-2}}]\label{lemma:characterization}
    Assume that $I:\mathcal{E}\to\mathbb{R}$ is of class $C^1$. Let $\{u_n\}$ be a sequence in $\mathcal{M}$ which is bounded in $\mathcal{E}$. Then the following are equivalent:
      \begin{itemize}
        \item[$(i)$] $\|I'_{|\mathcal{M}}(u_n)\|\to0$ as $n\to\infty$.
        \item[$(ii)$] $I'(u_n)-m^{-1}\langle I'(u_n),u_n\rangle u_n\to0$ in $\mathcal{E}^{-1}$ (the dual space of $\mathcal{E}$) as $n\to\infty$.
      \end{itemize}
    Here the last $u_n$ in $(ii)$ is an element of $\mathcal{E}^{-1}$ such that $\langle u_n,v\rangle:=(u_n,v)_\mathcal{H}$ for all $v\in\mathcal{E}$.
  \end{lemma}

\section{Preliminary results}\label{sect:preliminaries}

In this section we present some preliminary results.  For later convenience but without loss of generality, the exponent $q$ appeared in $(f3)$ will be denoted by $q_*$ and always understood as $2<q_*<2^*$.  The first technical result is Lemma \ref{lemma:geo1} below, which is a slightly modified version of \cite[Lemma 2.2]{Sh14}.
  \begin{lemma}\label{lemma:geo1}
    Assume that $N\geq1$ and $f$ satisfies $(f1)-(f3)$. Then the following statements hold.
      \begin{itemize}
        \item[$(i)$] Let $\{u_n\}$ be a bounded sequence in $H^1(\mathbb{R}^N)$. We have
            \begin{linenomath*}
              \begin{equation*}
                \underset{n\to\infty}{\lim}\int_{\mathbb{R}^N}F(u_n)dx=0
              \end{equation*}
            \end{linenomath*}
          if either $\lim_{n\to\infty}\|u_n\|_{L^2(\mathbb{R}^N)}=0$ or $\lim_{n\to\infty}\|u_n\|_{L^{q_*}(\mathbb{R}^N)}=0$.
        \item[$(ii)$] There exists $C=C(f,N,m)>0$ depending on $f$, $N$ and $m>0$ such that
             \begin{linenomath*}
               \begin{equation}\label{eq:geo1_1}
                 I(u)\geq\frac{1}{4}\int_{\mathbb{R}^N}|\nabla u|^2dx-C(f,N,m)
               \end{equation}
             \end{linenomath*}
           for all $u\in H^1(\mathbb{R}^N)$ satisfying $\|u\|^2_{L^2(\mathbb{R}^N)}\leq m$.
      \end{itemize}
  \end{lemma}

  \begin{lemma}\label{lemma:BL}
    Assume that $N\geq2$, $\{u_n\}\subset H^1(\mathbb{R}^N)$ is a bounded sequence, and $u_n\to u$ almost everywhere in $\mathbb{R}^N$ for some $u\in H^1(\mathbb{R}^N)$. Let $F:\mathbb{R}\to\mathbb{R}$ be a function of class $C^1$ with $F(0)=0$. If
      \begin{itemize}
        \item[$(i)$] when $N=2$, for any $\alpha>0$, there exists $C_\alpha>0$ such that
                       \begin{linenomath*}
                         \begin{equation*}
                           |F'(t)|\leq C_\alpha \left[|t|+\left(e^{\alpha t^2}-1\right)\right]\qquad\text{for all}~t\in\mathbb{R}.
                         \end{equation*}
                       \end{linenomath*}
        \item[$(ii)$] when $N\geq3$, there exists $C>0$ such that $|F'(t)|\leq C\left(|t|+|t|^{2^*-1}\right)$ for all $t\in\mathbb{R}$,
      \end{itemize}
    then
      \begin{linenomath*}
        \begin{equation}\label{eq:BL1}
          \lim_{n\to\infty}\int_{\mathbb{R}^N}\big|F(u_n)-F(u_n-u)-F(u)\big|dx=0.
        \end{equation}
      \end{linenomath*}
  \end{lemma}
\proof We prove this lemma by applying the Brezis-Lieb Lemma (\cite[Theorem 2]{BL83}). Clearly,
  \begin{linenomath*}
    \begin{equation}\label{eq:BL2}
      |F(a+b)-F(a)|=\left|\int^1_0\frac{d}{d\tau}F(a+\tau b)d\tau\right|=\left|\int^1_0F'(a+\tau b)bd\tau\right|\qquad\text{for all}~a,b\in\mathbb{R}.
    \end{equation}
  \end{linenomath*}
Let $\varepsilon>0$ be arbitrary. When $N=2$, by \eqref{eq:BL2}, $(i)$ and Young's inequality, one has
  \begin{linenomath*}
    \begin{equation*}
      \begin{split}
        |F(a+b)-F(a)|
        &\leq C_\alpha \int^1_0\left\{|a+\tau b|+\left[e^{\alpha \left(a+\tau b\right)^2}-1\right]\right\}|b|d\tau\\
        &\leq C_\alpha\left[|a|+|b|+\left(e^{4\alpha a^2}-1\right)+\left(e^{4\alpha b^2}-1\right)\right]|b|\\
        &\leq C_\alpha\left[\varepsilon a^2+\varepsilon^{-1}b^2+b^2+\varepsilon\left(e^{4\alpha a^2}-1\right)^2+\varepsilon^{-1} b^2+\left(e^{4\alpha b^2}-1\right)^2+b^2\right]\\
        &\leq \varepsilon C_\alpha\left[a^2+\left(e^{8\alpha a^2}-1\right)\right]+C_\alpha\left[2\left(1+\varepsilon^{-1}\right)b^2+\left(e^{8\alpha b^2}-1\right)\right]\\
        &=:\varepsilon\varphi(a)+\psi_\varepsilon(b).
      \end{split}
    \end{equation*}
  \end{linenomath*}
In particular, $|F(b)|\leq \psi_\varepsilon(b)$ for all $b\in\mathbb{R}$. Choose $M\geq1$ sufficiently large and $\alpha>0$ small enough such that
  \begin{linenomath*}
    \begin{equation*}
      \|u_n\|_{H^1(\mathbb{R}^2)},~\|u_n-u\|_{H^1(\mathbb{R}^2)},~\|u\|_{H^1(\mathbb{R}^2)}\leq M
    \end{equation*}
  \end{linenomath*}
and
  \begin{linenomath*}
    \begin{equation*}
      8\alpha\leq \frac{\beta}{M^2}\qquad\text{for some}~\beta\in(0,4\pi).
    \end{equation*}
  \end{linenomath*}
By the Moser-Trudinger inequality, we know that $\int_{\mathbb{R}^2}\varphi(u_n-u)dx$ is bounded uniformly in $\varepsilon$ and $n$, $\int_{\mathbb{R}^2}\psi_\varepsilon(u)dx<\infty$ for any $\varepsilon>0$, and $F(u)\in L^1(\mathbb{R}^2)$. In view of \cite[Theorem 2]{BL83}, we obtain \eqref{eq:BL1}.

When $N\geq3$, by \eqref{eq:BL2}, $(ii)$ and Young's inequality, one has
  \begin{linenomath*}
    \begin{equation*}
      \begin{split}
        |F(a+b)-F(a)|
        &\leq C\int^1_0\left(|a+\tau b|+|a+\tau b|^{2^*-1}\right)|b|d\tau\\
        &\leq C\left(|a|+2^{2^*}|a|^{2^*-1}+|b|+ 2^{2^*}|b|^{2^*-1}\right)|b|\\
        &\leq \varepsilon C\left(a^2+|2a|^{2^*}\right)+C\left[\left(1+\varepsilon^{-1}\right)b^2+\left(1+\varepsilon^{1-2^*}\right)|2b|^{2^*}\right]\\
        &=:\varepsilon\varphi(a)+\psi_\varepsilon(b).
      \end{split}
    \end{equation*}
  \end{linenomath*}
Applying Sobolev inequality and \cite[Theorem 2]{BL83}, one can conclude easily that \eqref{eq:BL1} holds.~~$\square$

  \begin{lemma}[{\cite[Corollary 3.2]{Me17}}]\label{lemma:lions}
    Assume that $N\geq4$ and $N-2M\neq0$. Let $\{u_n\}$ be a bounded sequence in $H^1_{\mathcal{O}_1}(\mathbb{R}^N)$ which satisfies
      \begin{linenomath*}
        \begin{equation}\label{eq:lions1}
          \lim_{r\to\infty}\left(\underset{n\to\infty}{\lim}\underset{y\in\{0\}\times\{0\}\times\mathbb{R}^{N-2M}}{\sup}\int_{B(y,r)}|u_n|^2dx\right)=0.
        \end{equation}
      \end{linenomath*}
    Then $u_n\to0$ in $L^p(\mathbb{R}^N)$ for any $2<p<2^*$.
  \end{lemma}
\proof We give here a complete proof for the reader's convenience. By \cite[Lemma 1.21]{Wi96}, the proof will be over if we can show that
  \begin{linenomath*}
    \begin{equation}\label{eq:lions2}
      \underset{n\to\infty}{\lim}\underset{y\in\mathbb{R}^N}{\sup}\int_{B(y,1)}|u_n|^2dx=0.
    \end{equation}
  \end{linenomath*}
We assume by contradiction that \eqref{eq:lions2} does not hold. Thus, up to a subsequence, there exist $\delta>0$ and $\{y_n\}\subset \mathbb{R}^N$ such that
  \begin{linenomath*}
    \begin{equation}\label{eq:lions3}
      \int_{B(y_n,1)}|u_n|^2dx\geq \delta>0\qquad\text{for}~n\geq1~\text{large enough}.
    \end{equation}
  \end{linenomath*}
	Since $\{u_n\}$ is bounded in $L^2(\mathbb{R}^N)$ and invariant with respect to $\mathcal{O}_1$, in view of \eqref{eq:lions3}, we deduce that $\{|(y^1_n,y^2_n)|\}$ must be bounded. Indeed, if $|(y^1_n,y^2_n)|\to\infty$, then one will derive  the existence of an arbitrarily large number of disjoint unit balls in the family $\{B(g^{-1}y_n,1)\}_{g\in\mathcal{O}_1}$.  Thus, for sufficiently large $r$, we have
  \begin{linenomath*}
    \begin{equation*}
      \int_{B((0,0,y^3_n),r)}|u_n|^2dx\geq\int_{B(y_n,1)}|u_n|^2dx\geq \delta>0,
    \end{equation*}
  \end{linenomath*}
which contradicts \eqref{eq:lions1}. Therefore, \eqref{eq:lions2} is satisfied and the desired conclusion follows.~~$\square$

For any $k\in\mathbb{N}$, let $\mathbb{S}^{k-1}$ be the unit sphere in $\mathbb{R}^k$, i.e.,
  \begin{linenomath*}
    \begin{equation*}
      \mathbb{S}^{k-1}:=\{\sigma\in\mathbb{R}^{k}~|~|\sigma|=1\}.
    \end{equation*}
  \end{linenomath*}
Recall that $X_2:=H^1_{\mathcal{O}_2}(\mathbb{R}^N)\cap X_\tau$. To proceed further, we need
  \begin{lemma}\label{lemma:keymapping}
    Assume that $N\geq 4$ and $f$ is an odd continuous function satisfying $(f4)$. Then, for any $k\in\mathbb{N}$, there exists an odd continuous mapping $\pi_k:\mathbb{S}^{k-1}\to X_2\setminus\{0\}$ such that
  \begin{linenomath*}
    \begin{equation*}
      \inf_{\sigma\in\mathbb{S}^{k-1}}\int_{\mathbb{R}^N}F(\pi_k[\sigma])dx\geq1\qquad\text{and}\qquad \sup_{\sigma\in\mathbb{S}^{k-1}}\|\pi_k[\sigma]\|_{L^\infty(\mathbb{R}^N)}\leq 2\zeta.
    \end{equation*}
  \end{linenomath*}
  \end{lemma}
\proof The first construction of such a mapping was done in \cite{JL18}, see \cite[Lemmas 4.2 and 4.3]{JL18}. For completeness we include here a shorter construction which borrows elements from \cite{Me17}.

Fix $k\in\mathbb{N}$. In view of \cite[Theorem 10]{Be83-2} and \cite[Proof of Lemma 8]{Be83-2}, there exist constants $R(k)>2(k+1)$ and $c_k>0$ such that, for any $R\geq R(k)$, there exists an odd continuous mapping $\tau_{k,R}:\mathbb{S}^{k-1} \to H^1(\mathbb{R}^N)$ having the properties that $\tau_{k,R}[\sigma]$ is a radial function, $\text{supp}\big(\tau_{k,R}[\sigma]\big)\subset \overline{B}(0,R)$ for any $\sigma \in \mathbb{S}^{k-1}$, $\sup_{\sigma\in\mathbb{S}^{k-1}}\|\tau_{k,R}[\sigma]\|_{L^\infty(\mathbb{R}^N)} = \zeta$ and
  \begin{linenomath*}
    \begin{equation}\label{eq:tau}
      \inf_{\sigma\in \mathbb{S}^{k-1}}\int_{\mathbb{R}^N}F\big(\tau_{k,R}[\sigma]\big)dx \geq c_k R^N.
    \end{equation}
  \end{linenomath*}
Let $\chi: \mathbb{R} \to [0,1]$ be an odd smooth function such that $\chi(t)=1$ for any $t \geq 1$. Following \cite[Remark 4.2]{Me17}, we define
  \begin{linenomath*}
    \begin{equation*}
      \pi_{k,R}[\sigma](x):=\tau_{k,R}[\sigma](x)\cdot \chi\big(|x_1|-|x_2|\big)
    \end{equation*}
  \end{linenomath*}
where $\sigma \in \mathbb{S}^{k-1}$ and $x=(x_1,x_2,x_3) \in \mathbb{R}^M \times \mathbb{R}^M \times \mathbb{R}^{N-2M}$. Here, we agree that the component $x_3$ does not exist when $N=2M$. Clearly, $\pi_{k,R}$ is an odd continuous mapping from $\mathbb{S}^{k-1}$ to $X_2$,
  \begin{linenomath*}
    \begin{equation*}
      \sup_{\sigma \in \mathbb{S}^{k-1}}\|\pi_{k,R}[\sigma]\|_{L^\infty(\mathbb{R}^N)} \leq \zeta \qquad \text{and}\qquad \text{supp}\big(\pi_{k,R}[\sigma]\big)\subset \overline{B}(0,R) \quad\text{for any}~\sigma \in \mathbb{S}^{k-1}.
    \end{equation*}
  \end{linenomath*}
Denoted by $\omega_l$ the surface area of $\mathbb{S}^l$. Set $A := \max_{t \in [0,\zeta]}|F(t)|$, $\omega:=\omega_{N-2M-1}\omega^2_{M-1}$ and $r_i:=|x_i|$ for $i=1,2,3$. For any $\sigma\in \mathbb{S}^{k-1}$, it is not difficult to see that
  \begin{linenomath}
    \begin{equation}\label{eq:pi}
      \begin{split}
        \int_{\mathbb{R}^N}F\big(\pi_{k,R}[\sigma]) dx
          &= \omega \int_0^{\infty} \int_0^{\infty} \int_0^{\infty} F\big(\pi_{k,R}[\sigma]) r_1^{M-1} r_2^{M-1} r_3^{N-2M-1} dr_1 dr_2 dr_3\\
          &= 2 \omega\int_0^{R} \int_0^{R} \int_{r_2}^{R} F\big(\pi_{k,R}[\sigma]) r_1^{M-1} r_2^{M-1} r_3^{N-2M-1} dr_1 dr_2 dr_3\\
		  &= 2 \omega\int_0^{R} \int_0^{R} \int_{r_2}^{R} F\big(\tau_{k,R}[\sigma]) r_1^{M-1} r_2^{M-1} r_3^{N-2M-1} dr_1 dr_2 dr_3\\
		  &\qquad -2 \omega\int_0^{R} \int_0^{R} \int_{r_2}^{r_2+1} F\big(\tau_{k,R}[\sigma]) r_1^{M-1} r_2^{M-1} r_3^{N-2M-1} dr_1 dr_2 dr_3\\
		  &\qquad + 2 \omega\int_0^{R} \int_0^{R} \int_{r_2}^{r_2+1} F\big(\pi_{k,R}[\sigma]\big) r_1^{M-1} r_2^{M-1} r_3^{N-2M-1} dr_1 dr_2 dr_3\\
          &\geq \int_{\mathbb{R}^N} F\big(\tau_{k,R}[\sigma]) dx - 2\omega A C_N R^{N-1},
      \end{split}
    \end{equation}
  \end{linenomath}
where $C_N >0$ is a constant depending only on $N$. Thus, in view of \eqref{eq:tau} and \eqref{eq:pi}, we obtain the desired mapping $\pi_k := \pi_{k,R}$ by taking $R> R(k)$ large enough.~~$\square$

  \begin{lemma}\label{lemma:geo2}
    Assume that $N\geq4$ and $f$ satisfies $(f1)-(f5)$. Then, for any $m>0$ and $k\in\mathbb{N}$, there exists an odd continuous mapping $\gamma_{m,k}:\mathbb{S}^{k-1}\to S_m\cap X_2$. Moreover, the following statements hold.
      \begin{itemize}
        \item[$(i)$] For any $k\in\mathbb{N}$, there exists  $m(k)> 0$ large enough  such that
                     \begin{linenomath*}
                       \begin{equation}\label{eq:geo2_1}
                         \sup_{\sigma\in\mathbb{S}^{k-1}}I(\gamma_{m,k}[\sigma])<0\qquad\text{for all}~m > m(k).
                       \end{equation}
                     \end{linenomath*}
				\item[$(ii)$] For any $s>0$, define $\gamma^s_{m,k}:\mathbb{S}^{k-1}\to S_m\cap X_2$ as follows:
                        \begin{linenomath*}
                          \begin{equation*}
                            \gamma^s_{m,k}[\sigma](x):=s^{N/2}\gamma_{m,k}[\sigma](sx),\qquad x\in\mathbb{R}^N~\text{and}~\sigma\in \mathbb{S}^{k-1}.
                          \end{equation*}
                        \end{linenomath*}
                      Then
                        \begin{linenomath*}
                          \begin{equation}\label{eq:geo2_2}
                            \limsup_{s\to0^+}\left(\sup_{\sigma\in\mathbb{S}^{k-1}}I(\gamma^s_{m,k}[\sigma])\right)\leq 0.
                          \end{equation}
                        \end{linenomath*}
                      If in addition \eqref{eq:f_key1} holds, then there exists $s_*>0$ small enough such that
                        \begin{linenomath*}
                          \begin{equation}\label{eq:geo2_3}
                            \sup_{\sigma\in\mathbb{S}^{k-1}}I(\gamma^s_{m,k}[\sigma])< 0\qquad\text{for any}~s\in(0,s_*).
                          \end{equation}
                        \end{linenomath*}
      \end{itemize}						
  \end{lemma}
\proof Fix $m>0$ and $k\in\mathbb{N}$. Using the mapping $\pi_k$ obtained in Lemma \ref{lemma:keymapping}, we can define an odd continuous mapping $\gamma_{m,k}:\mathbb{S}^{k-1}\to S_m\cap X_2$ as follows:
  \begin{linenomath*}
    \begin{equation*}
      \gamma_{m,k}[\sigma](x):=\pi_k[\sigma]\left(m^{-1/N}\cdot \|\pi_k[\sigma]\|^{2/N}_{L^2(\mathbb{R}^N)}\cdot x\right),\qquad x\in\mathbb{R}^N~\text{and}~\sigma\in\mathbb{S}^{k-1}.
    \end{equation*}
  \end{linenomath*}
Clearly,
  \begin{linenomath*}
    \begin{equation*}
      \sup_{\sigma\in \mathbb{S}^{k-1}}\|\gamma_{m,k}(\sigma)\|_{L^\infty(\mathbb{R}^N)}\leq 2\zeta.
    \end{equation*}
  \end{linenomath*}
We next show that this mapping $\gamma_{m,k}$ satisfies Items $(i)$ and $(ii)$.

$(i)$ Since $\mathbb{S}^{k-1}$ is compact and $0\notin \pi_k[\mathbb{S}^{k-1}]$, one can find $\alpha_k,\beta_k,\beta'_k>0$ independent of $\sigma\in\mathbb{S}^{k-1}$ such that
  \begin{linenomath*}
    \begin{equation*}
      \int_{\mathbb{R}^N}\bigl|\nabla\pi_k[\sigma]\bigr|^2dx\leq \alpha_k\qquad\text{and}\qquad \beta_k\leq\|\pi_k[\sigma]\|^2_{L^2(\mathbb{R}^N)}\leq \beta'_k.
    \end{equation*}
  \end{linenomath*}
Thus,
  \begin{linenomath*}
    \begin{equation*}
      \begin{split}
        I(\gamma_{m,k}[\sigma])&=\frac{1}{2}\int_{\mathbb{R}^N}\bigl|\nabla \gamma_{m,k}[\sigma]\bigr|^2dx-\int_{\mathbb{R}^N}F(\gamma_{m,k}[\sigma])dx\\
        &=\frac{m^{\frac{N-2}{N}}}{2\|\pi_k[\sigma]\|^{2(N-2)/N}_{L^2(\mathbb{R}^N)}}\int_{\mathbb{R}^N}\bigl|\nabla \pi_k[\sigma]\bigr|^2dx-\frac{m}{\|\pi_k[\sigma]\|^2_{L^2(\mathbb{R}^N)}}\int_{\mathbb{R}^N}F(\pi_k[\sigma])dx\\
        &\leq \frac{1}{2} \alpha_k \beta^{(2-N)/N}_k\cdot m^{\frac{N-2}{N}}-\big(\beta'_k\big)^{-1}\cdot m =:g_k(m).
      \end{split}
    \end{equation*}
  \end{linenomath*}
Clearly, $g_k(m)<0$ for sufficiently large $m>0$ and thus \eqref{eq:geo2_1} holds.

$(ii)$ We first prove \eqref{eq:geo2_2}. Let $\varepsilon>0$ be arbitrary. By $(f2)$, there exists $\delta>0$ such that
  \begin{linenomath*}
    \begin{equation*}
      |F(t)|\leq\varepsilon t^2\qquad\text{for all}~|t|\leq\delta.
    \end{equation*}
  \end{linenomath*}
Noting that
  \begin{linenomath*}
    \begin{equation}\label{eq:geo2_4}
      \|\gamma^s_{m,k}[\sigma]\|_{L^\infty(\mathbb{R}^N)}\leq 2s^{N/2}\zeta\qquad\text{for all}~\sigma\in\mathbb{S}^{k-1},
    \end{equation}
  \end{linenomath*}
one can find $s(\varepsilon)>0$ small enough such that
  \begin{linenomath*}
    \begin{equation*}
      \sup_{\sigma\in\mathbb{S}^{k-1}}\|\gamma^s_{m,k}[\sigma]\|_{L^\infty(\mathbb{R}^N)}\leq \delta\qquad\text{for all}~0<s<s(\varepsilon).
    \end{equation*}
  \end{linenomath*}
Therefore, for any $\sigma\in\mathbb{S}^{k-1}$ and $0<s<s(\varepsilon)$, we have
  \begin{linenomath*}
    \begin{equation*}
      \begin{split}
        I(\gamma^s_{m,k}[\sigma])&\leq \frac{1}{2}\int_{\mathbb{R}^N}\bigl|\nabla \gamma^s_{m,k}[\sigma]\bigr|^2dx+\int_{\mathbb{R}^N}\bigl|F(\gamma^s_{m,k}[\sigma])\bigr|dx\\
        &\leq\frac{1}{2}s^2\int_{\mathbb{R}^N}\bigl|\nabla \gamma_{m,k}[\sigma]\bigr|^2dx+\varepsilon \int_{\mathbb{R}^N} \bigl|\gamma_{m,k}[\sigma]\bigr|^2dx\\
        &=\frac{1}{2}s^2\int_{\mathbb{R}^N}\bigl|\nabla \gamma_{m,k}[\sigma]\bigr|^2dx+m\varepsilon.
       \end{split}
    \end{equation*}
  \end{linenomath*}
Since $\mathbb{S}^{k-1}$ is compact, there exists $C>0$ (independent of $\varepsilon>0$ and $s>0$) such that
  \begin{linenomath*}
    \begin{equation*}
      \sup_{\sigma\in\mathbb{S}^{k-1}}\int_{\mathbb{R}^N}\bigl|\nabla \gamma_{m,k}[\sigma]\bigr|^2dx\leq C.
    \end{equation*}
  \end{linenomath*}
Thus, for any $0<s<\min\bigl\{s(\varepsilon),(2\varepsilon/C)^{1/2}\bigr\}$, we obtain
  \begin{linenomath*}
    \begin{equation*}
      \sup_{\sigma\in\mathbb{S}^{k-1}}I(\gamma^s_{m,k}[\sigma])\leq \frac{1}{2}Cs^2+m\varepsilon\leq(m+1)\varepsilon.
    \end{equation*}
  \end{linenomath*}
Clearly, it follows that \eqref{eq:geo2_2} holds.

We now assume \eqref{eq:f_key1} and prove \eqref{eq:geo2_3}. For
  \begin{linenomath*}
    \begin{equation*}
      D:=\sup_{\sigma\in\mathbb{S}^{k-1}}\int_{\mathbb{R}^N}\bigl|\nabla \gamma_{m,k}[\sigma]\bigr|^2dx\biggm/\inf_{\sigma\in\mathbb{S}^{k-1}}\int_{\mathbb{R}^N}\bigl|\gamma_{m,k}[\sigma]\bigr|^{2+\frac{4}{N}}dx>0,
    \end{equation*}
  \end{linenomath*}
by \eqref{eq:f_key1}, there exists a $\delta>0$ such that
  \begin{linenomath*}
    \begin{equation*}
      F(t)\geq D|t|^{2+\frac{4}{N}}\qquad\text{for all}~|t|<\delta.
    \end{equation*}
  \end{linenomath*}
	Also, in  view of \eqref{eq:geo2_4}, one can find $s_*>0$ small enough  such that
  \begin{linenomath*}
    \begin{equation*}
      \sup_{\sigma\in\mathbb{S}^{k-1}}\|\gamma^s_{m,k}[\sigma]\|_{L^\infty(\mathbb{R}^N)}\leq \delta\qquad\text{for all}~0<s<s_*.
    \end{equation*}
  \end{linenomath*}
Thus, for any $\sigma\in\mathbb{S}^{k-1}$ and $0<s<s_*$, we have
  \begin{linenomath*}
    \begin{equation*}
      \begin{split}
        I(\gamma^s_{m,k}[\sigma])&\leq \frac{1}{2}\int_{\mathbb{R}^N}\bigl|\nabla \gamma^s_{m,k}[\sigma]\bigr|^2dx-D\int_{\mathbb{R}^N}\bigl|\gamma^s_{m,k}[\sigma]\bigr|^{2+\frac{4}{N}}dx\\
        &=\frac{1}{2}s^2\int_{\mathbb{R}^N}\bigl|\nabla \gamma_{m,k}[\sigma]\bigr|^2dx- Ds^2\int_{\mathbb{R}^N}\bigl|\gamma_{m,k}[\sigma]\bigr|^{2+\frac{4}{N}}dx\\
        &\leq-\frac{1}{2}s^2\sup_{\sigma\in\mathbb{S}^{k-1}}\int_{\mathbb{R}^N}\bigl|\nabla \gamma_{m,k}[\sigma]\bigr|^2dx<0,
       \end{split}
    \end{equation*}
  \end{linenomath*}
which implies \eqref{eq:geo2_3}.~~$\square$

  \begin{lemma}\label{lemma:PS}
    Assume that $N\geq4$, $N-2M\neq1$, and $f$ satisfies $(f1)-(f3)$ and $(f5)$. Then $I_{|S_m\cap X_2}$ satisfies the $(PS)_c$ condition for all $c<0$.
  \end{lemma}
\proof For given $c<0$, let $\{u_n\}\subset S_m\cap X_2$ be any sequence such that
  \begin{linenomath*}
    \begin{equation}\label{eq:PS1}
      I(u_n)\to c<0
    \end{equation}
  \end{linenomath*}
and
  \begin{linenomath*}
    \begin{equation}\label{eq:PS2}
      I'_{|S_m\cap X_2}(u_n)\to0.
    \end{equation}
  \end{linenomath*}
By \eqref{eq:PS1} and Lemma \ref{lemma:geo1} $(ii)$, we see that $\{u_n\}$ is bounded in $X_2$. Thus, up to a subsequence, we may assume that $u_n\rightharpoonup u$ in $X_2$ and $u_n\to u$ almost everywhere in $\mathbb{R}^N$ for some $u\in X_2$. In addition, thanks to \cite[Corollary 1.25]{Wi96}, $u_n\to u$ in $L^p(\mathbb{R}^N)$ for any $p\in(2, 2^*)$. Also, we know from \eqref{eq:PS2} and Lemma \ref{lemma:characterization} that
  \begin{linenomath*}
    \begin{equation}\label{eq:PS3}
      -\Delta u_n+\mu_n u_n-f(u_n)\to0\qquad\text{in}~(X_2)^{-1},
    \end{equation}
  \end{linenomath*}
where
  \begin{linenomath*}
    \begin{equation*}
      \mu_n:=\frac{1}{m}\left(\int_{\mathbb{R}^N}f(u_n)u_ndx-\int_{\mathbb{R}^N}|\nabla u_n|^2dx\right).
    \end{equation*}
  \end{linenomath*}
Since $\{\mu_n\}$ is bounded by $(f1)-(f3)$, we may assume that $\mu_n\to \mu$ for some $\mu\in\mathbb{R}$ and thus
  \begin{linenomath*}
    \begin{equation}
      -\Delta u+\mu u=f(u)\qquad\text{in}~(X_2)^{-1}.
    \end{equation}
  \end{linenomath*}
To show that $u_n\to u$ in $X_2$, the following two claims are needed.

\smallskip
\textbf{Claim 1.}
  \begin{linenomath*}
    \begin{equation}\label{eq:PS4}
      \lim_{n\to\infty}\int_{\mathbb{R}^N}f(u_n)u_ndx=\int_{\mathbb{R}^N}f(u)udx.
    \end{equation}
  \end{linenomath*}

Let $v_n:=u_n-u$. Clearly,
  \begin{linenomath*}
    \begin{equation*}
      \int_{\mathbb{R}^N}\left[f(u_n)u_n-f(u)u\right]dx=\int_{\mathbb{R}^N}f(u_n)v_ndx+\int_{\mathbb{R}^N}\left[f(u_n)-f(u)\right]udx.
    \end{equation*}
  \end{linenomath*}
Since $u_n\rightharpoonup u$ in $X_2$, one can show in a standard way that $\int_{\mathbb{R}^N}\big[f(u_n)-f(u)\big]udx\to0$. We estimate the remaining term $\int_{\mathbb{R}^N}f(u_n)v_ndx$. For any $\varepsilon>0$, there exists $C_\varepsilon>0$ such that
  \begin{linenomath*}
    \begin{equation*}
      |f(t)|\leq \varepsilon |t|+C_\varepsilon|t|^{q_*-1}\qquad\text{for all}~t\in\mathbb{R}.
    \end{equation*}
  \end{linenomath*}
Therefore, by H\"{o}lder inequality, we have
  \begin{linenomath*}
    \begin{equation*}
       \left|\int_{\mathbb{R}^N}f(u_n)v_ndx\right|\leq \varepsilon \|u_n\|_{L^2(\mathbb{R}^N)}\|v_n\|_{L^2(\mathbb{R}^N)}+C_\varepsilon\|u_n\|^{q_*-1}_{L^{q_*}(\mathbb{R}^N)}\|v_n\|_{L^{q_*}(\mathbb{R}^N)}.
    \end{equation*}
  \end{linenomath*}
Since $\|v_n\|_{L^{q_*}(\mathbb{R}^N)}\to0$ and $\varepsilon>0$ is arbitrary, we see that
  \begin{linenomath*}
    \begin{equation*}
      \lim_{n\to\infty}\int_{\mathbb{R}^N}f(u_n)v_ndx=0
    \end{equation*}
  \end{linenomath*}
and thus \eqref{eq:PS4} holds.

\smallskip
\textbf{Claim 2.} $\mu>0$.

Since $u_n-u\to0$ in $L^{q_*}(\mathbb{R}^N)$, we have $\int_{\mathbb{R}^N}F(u_n-u)dx\to0$ via Lemma \ref{lemma:geo1} $(i)$ and then $\int_{\mathbb{R}^N}F(u_n)dx\to\int_{\mathbb{R}^N}F(u)dx$ by Lemma \ref{lemma:BL}. In view of \eqref{eq:PS1}, we deduce that
  \begin{linenomath*}
    \begin{equation*}
      \begin{split}
        I(u)
        &=\frac{1}{2}\int_{\mathbb{R}^N}|\nabla u|^2dx-\int_{\mathbb{R}^N}F(u)dx\\
        &\leq \frac{1}{2}\lim_{n\to\infty}\int_{\mathbb{R}^N}|\nabla u_n|^2dx-\lim_{n\to\infty}\int_{\mathbb{R}^N}F(u_n)dx\\
        &=\lim_{n\to\infty}I(u_n)= c<0.
      \end{split}
    \end{equation*}
  \end{linenomath*}
Since, by Palais principle of symmetric criticality \cite{Pa79} and Poho\u{z}aev identity,
  \begin{linenomath*}
    \begin{equation*}
      P(u):=\frac{N-2}{2N}\int_{\mathbb{R}^N}|\nabla u|^2dx+\frac{1}{2}\mu\int_{\mathbb{R}^N}u^2dx-\int_{\mathbb{R}^N}F(u)dx=0,
    \end{equation*}
  \end{linenomath*}
we conclude that
  \begin{linenomath*}
    \begin{equation*}
      0>I(u)=I(u)-P(u)=\frac{1}{N}\int_{\mathbb{R}^N}|\nabla u|^2dx-\frac{1}{2}\mu\int_{\mathbb{R}^N}u^2dx.
    \end{equation*}
  \end{linenomath*}
Clearly, this implies that $\mu>0$.

With Claims 1 and 2 in hand, we can now show the strong convergence. By \eqref{eq:PS3}-\eqref{eq:PS4} and the fact that $\mu_n\to\mu>0$, we have
  \begin{linenomath*}
    \begin{equation*}
      \begin{split}
         \int_{\mathbb{R}^N}|\nabla u|^2dx+\mu\int_{\mathbb{R}^N}u^2dx
         &=\int_{\mathbb{R}^N}f(u)udx\\
         &=\lim_{n\to\infty}\int_{\mathbb{R}^N}f(u_n)u_ndx=\lim_{n\to\infty}\int_{\mathbb{R}^N}|\nabla u_n|^2dx+\mu m\\
         &\geq \int_{\mathbb{R}^N}|\nabla u|^2dx+\mu\int_{\mathbb{R}^N}u^2dx.
      \end{split}
    \end{equation*}
  \end{linenomath*}
Clearly, $\lim_{n\to\infty}\int_{\mathbb{R}^N}|\nabla u_n|^2dx=\int_{\mathbb{R}^N}|\nabla u|^2dx$, $\int_{\mathbb{R}^N}u^2dx=m$, and thus $u_n\to u$ in $X_2$.~~$\square$

\section{Proofs of the main results}\label{sect:proofs}
\subsection{Proof of Theorem \ref{theorem:nonradialsolutions}}\label{subsect:theorem1.2}
In this subsection, we shall use Theorem \ref{theorem:minimax} to prove Theorem \ref{theorem:nonradialsolutions}. Recall that $N\geq4$, $N-2M\neq1$ and $X_2:=H^1_{\mathcal{O}_2}(\mathbb{R}^N)\cap X_\tau$.  For any $m>0$ and $k\in\mathbb{N}$, we define
  \begin{linenomath*}
    \begin{equation*}
      A_{m,k}:=\gamma_{m,k}[\mathbb{S}^{k-1}],
    \end{equation*}
  \end{linenomath*}
where $\gamma_{m,k}$ is the odd continuous mapping given by Lemma \ref{lemma:geo2}. Clearly, $A_{m,k}$ is a closed symmetric set and $\mathcal{G}(A_{m,k})\geq k$ by Proposition \ref{proposition:genus} $(iii)$. Therefore, the class
  \begin{linenomath*}
    \begin{equation*}
      \Gamma_{m,k}:=\{A\in \Sigma(S_m\cap X_2)~|~\mathcal{G}(A)\geq k\}
    \end{equation*}
  \end{linenomath*}
is nonempty and the minimax value
\begin{linenomath*}
  \begin{equation*}
    E_{m,k}:=\inf_{A\in \Gamma_{m,k}}\sup_{u\in A}I(u)
  \end{equation*}
\end{linenomath*}
is well defined.
  \begin{lemma}\label{lemma:Emk}
    \begin{itemize}
	  \item[$(i)$] $-\infty<E_{m,k}\leq E_{m,k+1}\leq0$ for all $m>0$ and $k\in\mathbb{N}$.
      \item[$(ii)$] For any $k\in\mathbb{N}$, there exists $m(k) > 0$ large enough  such that $E_{m,k}<0$ for all $m > m(k)$.
      \item[$(iii)$] If in addition \eqref{eq:f_key1} holds, then $E_{m,k}<0$ for all $m>0$ and $k\in\mathbb{N}$.
	  \item[$(iv)$] For any $k\in\mathbb{N}$, the mapping $m\mapsto E_{m,k}$ is nonincreasing and continuous.
    \end{itemize}
  \end{lemma}
  \proof
	$(i)$ Since $\Gamma_{m,k+1}\subset\Gamma_{m,k}$ and $I$ is bounded from below on $S_m\cap X_2$ by Lemma \ref{lemma:geo1} $(ii)$, we have
  \begin{linenomath*}
    \begin{equation*}
      E_{m,k+1}\geq E_{m,k}>-\infty.
    \end{equation*}
  \end{linenomath*}
Let $\gamma^s_{m,k}$ be the odd continuous mapping given by Lemma \ref{lemma:geo2} $(ii)$. Clearly, $\gamma^s_{m,k}[\mathbb{S}^{k-1}]\in \Gamma_{m,k}$ and thus
  \begin{linenomath*}
    \begin{equation}\label{eq:Emk1}
      E_{m,k}\leq \sup_{\sigma\in\mathbb{S}^{k-1}}I(\gamma^s_{m,k}[\sigma])\qquad\text{for any}~s>0.
    \end{equation}
  \end{linenomath*}
In view of \eqref{eq:geo2_2}, we conclude that $E_{m,k}\leq 0$. The proof of Item $(i)$ is complete.

$(ii)$ This item follows from the fact that
  \begin{linenomath*}
    \begin{equation*}
      E_{m,k}\leq \sup_{\sigma\in\mathbb{S}^{k-1}}I(\gamma_{m,k}[\sigma])
    \end{equation*}
  \end{linenomath*}
and Lemma \ref{lemma:geo2} $(i)$.

$(iii)$ This item is a direct consequence of \eqref{eq:Emk1} and \eqref{eq:geo2_3}.

$(iv)$ Fix $k\in\mathbb{N}$. To prove the claim that the mapping $m\mapsto E_{m,k}$ is nonincreasing, we only need to show that, when $s>m>0$,
  \begin{linenomath*}
    \begin{equation}\label{eq:Emk2}
      E_{s,k}\leq E_{m,k}+\varepsilon \qquad\text{for any}~\varepsilon>0~\text{sufficiently small}.
    \end{equation}
  \end{linenomath*}
Clearly, \eqref{eq:Emk2} follows from Item $(i)$ if $E_{m,k}=0$. Thus, without loss of generality, we may assume that $E_{m,k}<0$. Let $\varepsilon\in(0,-E_{m,k}/2)$ be arbitrary. By the definition of $E_{m,k}$, there exists $A_{m,k}(\varepsilon)\in\Gamma_{m,k}$ such that
  \begin{linenomath*}
    \begin{equation}\label{eq:Emk3}
      \sup_{u\in A_{m,k}(\varepsilon)}I(u)\leq E_{m,k}+\varepsilon<0.
    \end{equation}
   \end{linenomath*}
Let
  \begin{linenomath*}
    \begin{equation*}
      B_{s,k}:=\{v(\cdot)=u(\cdot/t)~|~u\in A_{m,k}(\varepsilon)\},
    \end{equation*}
  \end{linenomath*}
where $t:=(s/m)^{1/N}>1$. It is clear that $B_{s,k}\in\Gamma_{s,k}$ and thus
  \begin{linenomath*}
    \begin{equation}\label{eq:Emk4}
      E_{s,k}\leq \sup_{v\in B_{s,k}}I(v)=\sup_{u\in A_{m,k}(\varepsilon)}I(u(\cdot/t)).
    \end{equation}
  \end{linenomath*}
For any $u\in A_{m,k}(\varepsilon)$, by \eqref{eq:Emk3} and the fact that $t>1$, we have
  \begin{linenomath*}
    \begin{equation*}
      I(u(\cdot/t))=t^NI(u)+\frac{1}{2}t^{N-2}(1-t^2)\int_{\mathbb{R}^N}|\nabla u|^2dx\leq I(u)\leq E_{m,k}+\varepsilon.
    \end{equation*}
 \end{linenomath*}
In view of \eqref{eq:Emk4}, we get the desired conclusion \eqref{eq:Emk2} and thus $E_{m,k}$ is nonincreasing in $m>0$.

We next show that the mapping $m\mapsto E_{m,k}$ is continuous. Let $s\in(0,m/2)$. Since $E_{m,k}$ is nonincreasing in $m>0$, we see that $E_{m-s,k}$ and $E_{m+s,k}$ are monotonic and bounded as $s\to0^+$. Therefore, $E_{m-s,k}$ and $E_{m+s,k}$ have limits as $s\to0^+$. Noting that $E_{m-s,k}\geq E_{m,k}\geq E_{m+s,k}$ for all $s\in(0,m/2)$, we have
  \begin{linenomath*}
    \begin{equation*}
      \lim_{s\to0^+}E_{m-s,k}\geq E_{m,k}\geq\lim_{s\to0^+}E_{m+s,k}.
    \end{equation*}
  \end{linenomath*}
To complete the proof, we only need to prove the reverse inequality, that is,
  \begin{linenomath*}
    \begin{equation*}
      \lim_{s\to0^+}E_{m-s,k}\leq E_{m,k}\leq\lim_{s\to0^+}E_{m+s,k}.
    \end{equation*}
  \end{linenomath*}

\smallskip
\textbf{Claim 1.} $\lim_{s\to0^+}E_{m-s,k}\leq E_{m,k}$.

Let $\varepsilon\in(0,1)$ be arbitrary. By the definition of $E_{m,k}$ and Item $(i)$, there exists $A_{m,k}(\varepsilon)\in\Gamma_{m,k}$ such that
  \begin{linenomath*}
    \begin{equation*}
      \sup_{u\in A_{m,k}(\varepsilon)}I(u)\leq E_{m,k}+\varepsilon\leq 1.
    \end{equation*}
  \end{linenomath*}
Since $A:=\cup_{0<\varepsilon<1}A_{m,k}(\varepsilon)\subset S_m\cap X_2$, we know from Lemma \ref{lemma:geo1} $(ii)$ that $A$ is a bounded set in $X_2$. We define
  \begin{linenomath*}
    \begin{equation*}
      B_{m-s,k}:=\{v(\cdot)=u(\cdot/t_s)~|~u\in A_{m,k}(\varepsilon)\},
    \end{equation*}
  \end{linenomath*}
where $t_s:=[(m-s)/m]^{1/N}>0$. Clearly, $B_{m-s,k}\in\Gamma_{m-s,k}$ and thus
  \begin{linenomath*}
    \begin{equation*}
      \begin{aligned}
        E_{m-s,k}\leq \sup_{v\in B_{m-s,k}}I(v)&=\sup_{u\in A_{m,k}(\varepsilon)}I(u(\cdot/t_s))\\
        &\leq \sup_{u\in A_{m,k}(\varepsilon)}I(u)+\sup_{u\in A_{m,k}(\varepsilon)}\left|I(u(\cdot/t_s))-I(u)\right|\\
        &\leq E_{m,k}+\varepsilon+\sup_{u\in A}\left|I(u(\cdot/t_s))-I(u)\right|.
      \end{aligned}
    \end{equation*}
  \end{linenomath*}
Since $\varepsilon\in(0,1)$ is arbitrary and $\sup_{u\in A}\left|I(u(\cdot/t_s))-I(u)\right|$ is independent of $\varepsilon\in(0,1)$, one will obtain Claim 1 if
  \begin{linenomath*}
    \begin{equation}\label{eq:Emk5}
      \lim_{s\to0^+}\left(\sup_{u\in A}\left|I(u(\cdot/t_s))-I(u)\right|\right)=0.
    \end{equation}
  \end{linenomath*}
We now prove \eqref{eq:Emk5}. Noting that $t_s$ is only dependent on $s$, we have
  \begin{linenomath*}
    \begin{equation*}
      \sup_{u\in A}\left|I(u(\cdot/t_s))-I(u)\right|\leq
       \frac{1}{2}\left|t^{N-2}_s-1\right|\sup_{u\in A}\int_{\mathbb{R}^N}|\nabla u|^2dx+\left|t^N_s-1\right|\sup_{u\in A}\int_{\mathbb{R}^N}|F(u)|dx.
    \end{equation*}
  \end{linenomath*}
Since $A$ is bounded in $X_2$, by $(f1)-(f3)$, we see that $\sup_{u\in A}\int_{\mathbb{R}^N}|F(u)|dx$ is bounded uniformly in $s\in(0,m/2)$. In view of the fact that $\lim_{s\to0^+}t_s=1$, we obtain \eqref{eq:Emk5}.

\smallskip
\textbf{Claim 2.} $\lim_{s\to0^+}E_{m+s,k}\geq E_{m,k}$.

For any $s\in(0,m/2)$, by the definition of $E_{m+s,k}$, there exists $A_{m+s,k}\in\Gamma_{m+s,k}$ such that
  \begin{linenomath*}
    \begin{equation*}
      \sup_{u\in A_{m+s,k}}I(u)\leq E_{m+s,k}+s.
    \end{equation*}
  \end{linenomath*}
Let $A:=\cup_{0<s<m/2}A_{m+s,k}$. Since $\sup_{u\in A}\|u\|^2_{L^2(\mathbb{R}^N)}\leq 3m/2$ and $\sup_{u\in A}I(u)\leq m/2$ by Item $(i)$, we know from Lemma \ref{lemma:geo1} $(ii)$ that $A$ is a bounded set in $X_2$. Define
  \begin{linenomath*}
    \begin{equation*}
      B_{m,k}(s):=\{v(\cdot)=u(\cdot/t_s)~|~u\in A_{m+s,k}\},
    \end{equation*}
  \end{linenomath*}
where $t_s:=[m/(m+s)]^{1/N}>0$. Clearly, $B_{m,k}(s)\in\Gamma_{m,k}$ and thus
  \begin{linenomath*}
    \begin{equation*}
      \begin{aligned}
        E_{m,k}\leq \sup_{v\in B_{m,k}(s)}I(v)&=\sup_{u\in A_{m+s,k}}I(u(\cdot/t_s))\\
        &\leq \sup_{u\in A_{m+s,k}}I(u)+\sup_{u\in A_{m+s,k}}\left|I(u(\cdot/t_s))-I(u)\right|\\
        &\leq E_{m+s,k}+s+\sup_{u\in A}\left|I(u(\cdot/t_s))-I(u)\right|.
      \end{aligned}
    \end{equation*}
  \end{linenomath*}
Arguing as the proof of \eqref{eq:Emk5}, we have that
  \begin{linenomath*}
    \begin{equation*}
      \lim_{s\to0^+}\left(\sup_{u\in A}\left|I(u(\cdot/t_s))-I(u)\right|\right)=0.
    \end{equation*}
  \end{linenomath*}
Therefore,
  \begin{linenomath*}
    \begin{equation*}
      E_{m,k}\leq \lim_{s\to0^+}\left(E_{m+s,k}+s+\sup_{u\in A}\left|I(u(\cdot/t_s))-I(u)\right|\right)=\lim_{s\to0^+}E_{m+s,k}.
    \end{equation*}
  \end{linenomath*}
The proof of Claim 2 is complete.
~~$\square$

\medskip
\noindent
\textbf{Proof of Theorem \ref{theorem:nonradialsolutions}.}
Clearly, Theorem \ref{theorem:nonradialsolutions} $(i)$ and $(ii)$ are Lemma \ref{lemma:Emk} $(i)$ and $(iv)$ respectively. Let $\mathcal{E}:=X_2$ and $\mathcal{H}:=L^2(\mathbb{R}^N)$. For any $k\in\mathbb{N}$, we define
  \begin{linenomath*}
    \begin{equation*}
      m_k:=\inf\{m>0~|~E_{m,k}<0\}.
    \end{equation*}
  \end{linenomath*}
By Lemma \ref{lemma:Emk} $(i)$, $(ii)$ and $(iv)$, it follows that $m_k \in [0,\infty)$,
  \begin{linenomath*}
    \begin{equation*}
      E_{m,k}=0\quad\text{if}~0<m\leq m_k,\qquad E_{m,k}<0\quad\text{when}~m>m_k.
    \end{equation*}
  \end{linenomath*}
Fixing $k\in\mathbb{N}$, when $m>m_k$, we have
  \begin{linenomath*}
    \begin{equation*}
      -\infty<E_{m,1}\leq E_{m,2}\leq \cdots\leq E_{m,k}<0.
    \end{equation*}
  \end{linenomath*}
In view of Lemma \ref{lemma:geo1} $(ii)$, Lemma \ref{lemma:PS}, and Theorem \ref{theorem:minimax} $(i)$ and $(ii)$, we know that $I_{|S_m\cap X_2}$ has at least $k$ distinct critical points associated to the levels $E_{m,j}$ ($j= 1, 2, \cdots, k$). Thus, by Palais principle of symmetric criticality, we obtain Theorem \ref{theorem:nonradialsolutions} $(iii)$. If \eqref{eq:f_key1} holds, by Lemma \ref{lemma:Emk} $(iii)$, we see that $E_{m,k}<0$ for any $m>0$ and $k \in \mathbb{N}$ (and thus $m_k=0$ for any $k\in\mathbb{N}$). Applying Theorem \ref{theorem:minimax} $(i)$ and $(iii)$ to $I_{|S_m\cap X_2}$, we get Theorem \ref{theorem:nonradialsolutions} $(iv)$. ~~$\square$

\subsection{Proof of Theorem \ref{theorem:nonradialsolution}}\label{subsect:theorem1.1}
This subsection is devoted to the proof of Theorem \ref{theorem:nonradialsolution}.  Recall that $N\geq4$, $X_1:=H^1_{\mathcal{O}_1}(\mathbb{R}^N)\cap X_\tau$ and
  \begin{linenomath*}
    \begin{equation*}
      E_m:=\inf_{u\in S_m\cap X_1}I(u).
    \end{equation*}
  \end{linenomath*}
	Clearly, $E_m>-\infty$ by Lemma \ref{lemma:geo1} $(ii)$. Since $X_2\subset X_1$, using Lemma \ref{lemma:geo2} and arguing as the proof of Lemma \ref{lemma:Emk}, we also have
	\begin{lemma}\label{lemma:E1m}
    \begin{itemize}
      \item[$(i)$] $-\infty<E_m\leq0$ for all $m>0$.
      \item[$(ii)$] There exists $m_0>0$ large enough such that $E_m<0$ for all $m > m_0$.
      \item[$(iii)$] If in addition \eqref{eq:f_key1} holds, then $E_m<0$ for all $m>0$.
			\item[$(iv)$] The mapping $m\mapsto E_m$ is nonincreasing and continuous.
    \end{itemize}
  \end{lemma}
	
\begin{lemma}\label{lemma:E1m_1}
  For any $m>s>0$, one has
    \begin{linenomath*}
      \begin{equation}\label{eq:E1m3}
        E_m\leq \frac{m}{s}E_s.
      \end{equation}
    \end{linenomath*}
  If $E_s$ is reached, then the inequality is strict.
\end{lemma}
\proof  Let $t:=m/s>1$. For any $\varepsilon>0$, there exists $u\in S_s\cap X_1$ such that
  \begin{linenomath*}
    \begin{equation*}
      I(u)\leq E_s+\varepsilon.
    \end{equation*}
  \end{linenomath*}
Clearly, $w:=u(t^{-1/N}\cdot)\in S_m\cap X_1$ and then
  \begin{linenomath*}
    \begin{equation}\label{eq:E1m4}
      E_m\leq I(w)=tI(u)+\frac{1}{2}t^{\frac{N-2}{N}}\left(1-t^{\frac{2}{N}}\right)\int_{\mathbb{R}^N}|\nabla u|^2dx< tI(u)\leq\frac{m}{s}(E_s+\varepsilon).
    \end{equation}
  \end{linenomath*}
Since $\varepsilon>0$ is arbitrary, we see that \eqref{eq:E1m3} holds. If $E_s$ is reached, for example, at some $u\in S_s\cap X_1$, then we can let $\varepsilon=0$ in \eqref{eq:E1m4} and thus the strict inequality follows.~~$\square$

\medskip
\noindent
\textbf{Proof of Theorem \ref{theorem:nonradialsolution}.}
 We define
  \begin{linenomath*}
    \begin{equation}\label{eq:m^*}
      m^*:=\inf\{m>0~|~E_m<0\}.
    \end{equation}
  \end{linenomath*}
By Lemma \ref{lemma:E1m}, it is clear that $m^*\in[0,\infty)$,
  \begin{linenomath*}
    \begin{equation}\label{eq:E1m1}
      E_m=0\quad\text{if}~0<m\leq m^*,\qquad E_m<0\quad\text{when}~m>m^*;
    \end{equation}
  \end{linenomath*}
in particular, $m^*=0$ if \eqref{eq:f_key1} holds. Let us show that  if $0<m< m^*$, then $E_m=0$ is not reached. Indeed, assuming by contradiction that $E_m$ is reached for some $m\in(0,m^*)$, in view of Lemma \ref{lemma:E1m_1}, we have
  \begin{linenomath*}
    \begin{equation*}
      E_{m^*}<\frac{m^*}{m}E_m=0
    \end{equation*}
  \end{linenomath*}
which leads to a contradiction since $E_{m^*}=0$.

To complete the proof of Theorem \ref{theorem:nonradialsolution}, the only remaining task is to show that the infimum $E_m$ is reached when $m>m^*$.  When $N-2M=0$, we have $X_1=X_2$ (with $N-2M \neq 1$). Since in that case $E_m = E_{m,1}$ and $m^*= m_1$, the result follows directly from the property, established in Theorem \ref{theorem:nonradialsolutions}, that $E_{m,1}$ is a critical value. The rest of the proof is devoted to deal with the delicate case, that is when $N-2M \neq 0$.

Fix $m>m^*$ and let $\{u_n\}\subset S_m\cap X_1$ be a minimizing sequence with respect to $E_m$. Clearly, $\{u_n\}$ is bounded in $X_1$ by Lemma \ref{lemma:geo1} $(ii)$. Up to a subsequence, we may assume that $\lim_{n\to\infty}\int_{\mathbb{R}^N}|\nabla u_n|^2dx$ and $\lim_{n\to\infty}\int_{\mathbb{R}^N}F(u_n)dx$ exist. Since $E_m<0$ by \eqref{eq:E1m1}, we have that
  \begin{linenomath*}
    \begin{equation}\label{eq:nonvanishing_1}
      \lim_{r\to\infty}\left(\underset{n\to\infty}{\lim}\underset{y\in\{0\}\times\{0\}\times\mathbb{R}^{N-2M}}{\sup}\int_{B(y,r)}|u_n|^2dx\right)>0.
    \end{equation}
  \end{linenomath*}
Indeed, if \eqref{eq:nonvanishing_1} does not hold, then $u_n\to0$ in $L^{q_*}(\mathbb{R}^N)$ by Lemma \ref{lemma:lions} and thus
  \begin{linenomath*}
    \begin{equation*}
      \lim_{n\to\infty}\int_{\mathbb{R}^N}F(u_n)dx=0
    \end{equation*}
  \end{linenomath*}
via Lemma \ref{lemma:geo1} $(i)$; since $I(u_n)\geq -\int_{\mathbb{R}^N}F(u_n)dx$, a contradiction is obtained as follows:
  \begin{linenomath*}
    \begin{equation*}
      0>E_m=\lim_{n\to\infty}I(u_n)\geq -\lim_{n\to\infty}\int_{\mathbb{R}^N}F(u_n)dx=0.
    \end{equation*}
  \end{linenomath*}
With \eqref{eq:nonvanishing_1} in hand, we see that there exist $r_0>0$ and $\{y_n\}\subset\{0\}\times\{0\}\times\mathbb{R}^{N-2M}$ such that
  \begin{linenomath*}
    \begin{equation}\label{eq:nonvanishing_2}
      \underset{n\to\infty}{\lim}\int_{B(y_n,r_0)}|u_n|^2dx>0.
    \end{equation}
  \end{linenomath*}
Since $\{u_n(\cdot-y_n)\}\subset S_m\cap X_1$ is bounded, up to a subsequence, we may assume that $u_n(\cdot-y_n)\rightharpoonup u$ in $X_1$ for some $u\in X_1$, $u_n(\cdot-y_n)\to u$ in $L^2_{\text{loc}}(\mathbb{R}^N)$ and $u_n(\cdot-y_n)\to u$ almost everywhere in $\mathbb{R}^N$. Clearly, $u\neq0$ by \eqref{eq:nonvanishing_2}. Let $v_n:=u_n(\cdot-y_n)-u$. Noting that $v_n\rightharpoonup0$ in $X_1$, we have
  \begin{linenomath*}
    \begin{equation*}
      \int_{\mathbb{R}^N}|u+v_n|^2dx=\int_{\mathbb{R}^N}|u|^2dx+\int_{\mathbb{R}^N}|v_n|^2dx+o_n(1)
    \end{equation*}
  \end{linenomath*}
and
  \begin{linenomath*}
    \begin{equation*}
      \int_{\mathbb{R}^N}|\nabla (u+v_n)|^2dx=\int_{\mathbb{R}^N}|\nabla u|^2dx+\int_{\mathbb{R}^N}|\nabla v_n|^2dx+o_n(1),
    \end{equation*}
  \end{linenomath*}
where $o_n(1)\to0$ as $n\to\infty$. With the aid of Lemma \ref{lemma:BL}, we also have
  \begin{linenomath*}
    \begin{equation*}
       \lim_{n\to\infty}\int_{\mathbb{R}^N}F(u+v_n)dx=\int_{\mathbb{R}^N}F(u)dx+\lim_{n\to\infty}\int_{\mathbb{R}^N}F(v_n)dx.
    \end{equation*}
  \end{linenomath*}
Since $I(u_n)=I(u_n(\cdot-y_n))=I(u+v_n)$, it follows that
  \begin{linenomath*}
    \begin{equation}\label{eq:key_1}
      m=\|u\|^2_{L^2(\mathbb{R}^N)}+\lim_{n\to\infty}\|v_n\|^2_{L^2(\mathbb{R}^N)}
    \end{equation}
  \end{linenomath*}
and
  \begin{linenomath*}
    \begin{equation}\label{eq:key_2}
      E_m=\lim_{n\to\infty}I(u+v_n)=I(u)+\lim_{n\to\infty}I(v_n).
    \end{equation}
  \end{linenomath*}

We prove below a claim and then conclude the proof.

\smallskip
\textbf{Claim.}  $\lim_{n\to\infty}\|v_n\|^2_{L^2(\mathbb{R}^N)}=0$. In particular, by \eqref{eq:key_1}, $\|u\|^2_{L^2(\mathbb{R}^N)}=m$.

Let $t_n:=\|v_n\|^2_{L^2(\mathbb{R}^N)}$ for every $n\in\mathbb{N}$. If we assume that $\lim_{n\to\infty}t_n>0$, then \eqref{eq:key_1} implies that  $s:=\|u\|^2_{L^2(\mathbb{R}^N)}\in (0, m)$.  By the definition of $E_{t_n}$ and Lemma \ref{lemma:E1m} $(iv)$, we have
  \begin{linenomath*}
    \begin{equation*}
      \lim_{n\to\infty}I(v_n)\geq \lim_{n\to\infty}E_{t_n}= E_{m-s}.
    \end{equation*}
  \end{linenomath*}
From \eqref{eq:key_2} and Lemma \ref{lemma:E1m_1}, it follows that
\begin{linenomath*}
    \begin{equation*}
      E_m\geq I(u) + E_{m-s} \geq E_s + E_{m-s} \geq \frac{s}{m}E_m + \frac{m-s}{m}E_m= E_m.
    \end{equation*}
  \end{linenomath*}
Thus necessarily $I(u) = E_s$ and this shows that $E_s$ is reached at $u$. But then still from \eqref{eq:key_2} and Lemma \ref{lemma:E1m_1}, one has
\begin{linenomath*}
    \begin{equation*}
      E_m \geq E_s + E_{m-s} > \frac{s}{m}E_m + \frac{m-s}{m}E_m = E_m
    \end{equation*}
  \end{linenomath*}
which is a contradiction and thus proves the Claim.

\smallskip
\textbf{Conclusion.} Clearly, $u\in S_m\cap X_1$ by the Claim  and thus $I(u)\geq E_m$. Now since the Claim and Lemma \ref{lemma:geo1} $(i)$ imply that
  \begin{linenomath*}
    \begin{equation*}
      \lim_{n\to\infty}\int_{\mathbb{R}^N}F(v_n)dx=0,
    \end{equation*}
  \end{linenomath*}
we also have $\lim_{n \to \infty} I(v_n) \geq 0$. Thus, by \eqref{eq:key_2}, we get $E_m \geq I(u)$, hence $E_m$ is reached at $u\in S_{m}\cap X_1$. This completes the proof of Theorem \ref{theorem:nonradialsolution}.~~$\square$

  \begin{remark}\label{remark:positive}
    In the proof of Theorem \ref{theorem:nonradialsolution}  we define the number $m^*$ via \eqref{eq:m^*}. When we do not have \eqref{eq:f_key1}, this number can be positive.  To see this, following closely \cite{Sh14}, we assume, in addition to $(f1)-(f5)$, that
      \begin{linenomath*}
        \begin{equation}\label{eq:f_key2}
          \limsup_{t\to0}\frac{F(t)}{|t|^{2+\frac{4}{N}}}<+\infty.
        \end{equation}
      \end{linenomath*}
    Since there exists $C(f)>0$ such that $F(t)\leq C(f)|t|^{2+4/N}$ for any $t\in\mathbb{R}$, by Gagliardo-Nirenberg inequality, it follows that
      \begin{linenomath*}
        \begin{equation*}
          \int_{\mathbb{R}^N}F(u)dx\leq C(f)\int_{\mathbb{R}^N}|u|^{2+4/N}dx\leq C(f)C(N)m^{2/N}\int_{\mathbb{R}^N}|\nabla u|^2dx
        \end{equation*}
      \end{linenomath*}
    for all $u\in S_m$. Then, for any $m>0$ small enough such that $C(f)C(N)m^{2/N}\leq 1/4$, we have
       \begin{linenomath*}
         \begin{equation*}
           I(u)\geq\frac{1}{2}\int_{\mathbb{R}^N}|\nabla u|^2dx-\frac{1}{4}\int_{\mathbb{R}^N}|\nabla u|^2dx=\frac{1}{4}\int_{\mathbb{R}^N}|\nabla u|^2dx>0.
         \end{equation*}
       \end{linenomath*}
    Clearly, this implies that $E_m\geq0$ when $m>0$ is small.
  \end{remark}

\section{Multiple radial solutions}\label{sect:theoremB}
Based on the approach developed to prove Theorem \ref{theorem:nonradialsolutions} and under weak conditions, we give in this section a new proof for the result due to Hirata and Tanaka \cite{HT18} on multiple radial solutions.

\begin{theorem}\label{theorem:radialsolutions}
    Assume that $N\geq2$ and $f$ satisfies $(f1),(f2),(f3)',(f4)$ and $(f5)$.  Then the following statements hold.
     \begin{itemize}
       \item[$(i)$] For each $k\in\mathbb{N}$ there exists $\overline{m}_k\in[0,\infty)$ such that, when $m>\overline{m}_k$, \eqref{problem} has at least $k$ radial solutions (with negative energies).
       \item[$(ii)$] Assume in addition \eqref{eq:f_key1}, then \eqref{problem} has infinitely many radial solutions $\{v_n\}^\infty_{n=1}$ for all $m>0$. In particular, $I(v_n)<0$ for each $n\in\mathbb{N}$ and $I(v_n)\to0$ as $n\to\infty$.
     \end{itemize}
  \end{theorem}
Note that in \cite[Theorem 0.2]{HT18},  instead of $(f3)'$, it is required the stronger condition
\begin{itemize}
    \item[$(f3)''$] $\lim_{t\rightarrow\infty}f(t)/|t|^{1+4/N}=0$.
  \end{itemize}
To derive their result, Hirata and Tanaka apply a version of \emph{symmetric mountain pass argument} to $I(\lambda,u):\mathbb{R}\times H^1_r(\mathbb{R}^N)\to\mathbb{R}$, a Lagrange formulation of \eqref{problem} defined as
  \begin{linenomath*}
    \begin{equation*}
      I(\lambda,u)=\frac{1}{2}\int_{\mathbb{R}^N}|\nabla u|^2dx-\int_{\mathbb{R}^N}F(u)dx+\frac{1}{2}e^\lambda\left(\int_{\mathbb{R}^N}u^2dx-m\right).
    \end{equation*}
  \end{linenomath*}
Here $H^1_r(\mathbb{R}^N)$ stands for the space of radially symmetric functions in $H^1(\mathbb{R}^N)$.
Note also that  just assuming $(f1),(f2),(f3)''$ and $(f4)$, they established  the existence of one radial solution via a \emph{mountain pass argument} applied to $I(\lambda,u)$. As a consequence they derived a minimax characterization of the global infimum $E_m$, see \cite[Theorem 0.1]{HT18} for more details.

To prove Theorem \ref{theorem:radialsolutions}, in view of Remark \ref{remark:extension}, we can assume without loss of generality that $(f3)$ holds. Since when $N\geq2$ the embedding $H^1_r(\mathbb{R}^N)\hookrightarrow L^p(\mathbb{R}^N)$ is compact for any $2<p <2^*$, modifying the proof of Lemma \ref{lemma:PS} accordingly, we have the following compactness result.
  \begin{lemma}\label{lemma:PS1}
    The constrained functional $I_{|S_m\cap H^1_r(\mathbb{R}^N)}$ satisfies the $(PS)_c$ condition for all $c<0$.
  \end{lemma}

Since the remaining arguments are similar to that for Theorem \ref{theorem:nonradialsolutions}, we just outline the proof.  Fix $m>0$ and $k\in\mathbb{N}$. By \cite[Theorem 10]{Be83-2}, there exists an odd continuous mapping $\overline{\pi}_k:\mathbb{S}^{k-1}\to H^1_r(\mathbb{R}^N)\setminus\{0\}$ such that
  \begin{linenomath*}
    \begin{equation*}
      \inf_{\sigma\in\mathbb{S}^{k-1}}\int_{\mathbb{R}^N}F(\overline{\pi}_k[\sigma])dx\geq1\qquad\text{and}\qquad \sup_{\sigma\in\mathbb{S}^{k-1}}\|\overline{\pi}_k[\sigma]\|_{L^\infty(\mathbb{R}^N)}\leq \zeta.
    \end{equation*}
  \end{linenomath*}
Thus, an odd continuous mapping $\overline{\gamma}_{m,k}:\mathbb{S}^{k-1}\to S_m\cap H^1_r(\mathbb{R}^N)$ can be introduced as follows:
  \begin{linenomath*}
    \begin{equation*}
      \overline{\gamma}_{m,k}[\sigma](x):=\overline{\pi}_k[\sigma]\left(m^{-1/N}\cdot \|\overline{\pi}_k[\sigma]\|^{2/N}_{L^2(\mathbb{R}^N)}\cdot x\right),\qquad x\in\mathbb{R}^N~\text{and}~\sigma\in\mathbb{S}^{k-1}.
    \end{equation*}
  \end{linenomath*}
For any $s>0$, we then define $\overline{\gamma}^s_{m,k}[\sigma](x):=s^{N/2}\overline{\gamma}_{m,k}[\sigma](sx)$. Arguing as the proof of Lemma \ref{lemma:geo2}, we see that $\overline{\gamma}_{m,k}$ and $\overline{\gamma}^s_{m,k}$ satisfy the following lemma.
  \begin{lemma}\label{lemma:geo3}
    \begin{itemize}
      \item[$(i)$] For any $k\in\mathbb{N}$, there exists $\overline{m}(k) > 0$ large enough  such that
                     \begin{linenomath*}
                       \begin{equation*}
                         \sup_{\sigma\in\mathbb{S}^{k-1}}I(\overline{\gamma}_{m,k}[\sigma])<0\qquad\text{for all}~m > \overline{m}(k).
                       \end{equation*}
                     \end{linenomath*}
      \item[$(ii)$] For any $m>0$ and $k\in\mathbb{N}$, we have
                      \begin{linenomath*}
                        \begin{equation*}
                          \limsup_{s\to0^+}\left(\sup_{\sigma\in\mathbb{S}^{k-1}}I(\overline{\gamma}^s_{m,k}[\sigma])\right)\leq 0.
                        \end{equation*}
                      \end{linenomath*}
                    If in addition \eqref{eq:f_key1} holds, then there exists $s^*>0$ small enough such that
                      \begin{linenomath*}
                        \begin{equation*}
                          \sup_{\sigma\in\mathbb{S}^{k-1}}I(\overline{\gamma}^s_{m,k}[\sigma])< 0\qquad\text{for any}~s\in(0,s^*).
                        \end{equation*}
                      \end{linenomath*}
      \end{itemize}
  \end{lemma}

Since $\mathcal{G}(\overline{\gamma}_{m,k}[\mathbb{S}^{k-1}])\geq k$,  the class
  \begin{linenomath*}
    \begin{equation*}
      \overline{\Gamma}_{m,k}:=\{A\in \Sigma(S_m\cap H^1_r(\mathbb{R}^N))~|~\mathcal{G}(A)\geq k\}
    \end{equation*}
  \end{linenomath*}
is nonempty and the minimax value
  \begin{linenomath*}
    \begin{equation*}
      \overline{E}_{m,k}:=\inf_{A\in \overline{\Gamma}_{m,k}}\sup_{u\in A}I(u)
    \end{equation*}
  \end{linenomath*}
is well defined. With the aid of Lemma \ref{lemma:geo3}, repeating the argument of Lemma \ref{lemma:Emk}, we obtain
  \begin{lemma}\label{lemma:Ebarmk}
    \begin{itemize}
      \item[$(i)$] $-\infty<\overline{E}_{m,k}\leq \overline{E}_{m,k+1}\leq0$ for all $m>0$ and $k\in\mathbb{N}$.
      \item[$(ii)$] For any $k\in\mathbb{N}$, there exists $\overline{m}(k) > 0$ large enough  such that $\overline{E}_{m,k}<0$ for all $m > \overline{m}(k)$.
			\item[$(iii)$] If in addition \eqref{eq:f_key1} holds, then $\overline{E}_{m,k}<0$ for all $m>0$ and $k\in\mathbb{N}$.
      \item[$(iv)$] For any $k\in\mathbb{N}$, the mapping $m\mapsto \overline{E}_{m,k}$ is nonincreasing and continuous.
      \end{itemize}
  \end{lemma}

\medskip
\noindent
\textbf{Conclusion.} Let $\mathcal{E}:=H^1_r(\mathbb{R}^N)$ and $\mathcal{H}:=L^2(\mathbb{R}^N)$. For any $k\in\mathbb{N}$, define
  \begin{linenomath*}
    \begin{equation*}
      \overline{m}_k:=\inf\{m>0~|~\overline{E}_{m,k}<0\}.
    \end{equation*}
  \end{linenomath*}
In view of Lemma \ref{lemma:geo1} $(ii)$, Lemma \ref{lemma:PS1}, Lemma \ref{lemma:Ebarmk} and Theorem \ref{theorem:minimax}, we obtain Theorem \ref{theorem:radialsolutions}.~~$\square$



\end{document}